\documentclass[AMA,STIX1COL]{WileyNJD-v2}

\usepackage{multirow}
\usepackage{float}
\usepackage{makecell}
\usepackage{graphicx}
\usepackage{mathtools}
\usepackage{url}
\usepackage{textcomp,gensymb}
\usepackage{algorithm,algpseudocode} 
\usepackage{xcolor}
\usepackage{booktabs}
\usepackage{amsmath}
\usepackage{caption}
\usepackage{subcaption}
\newcommand{\tc}[2]{\textcolor{#1}{#2}}
\newcommand{\change}[1]{\tc{black}{#1}}

\colorlet{blue}{blue!80!green}
\colorlet{red}{red!80!green}

\global\long\def\norm#1{\lVert#1\rVert}%
\global\long\def\abs#1{\lvert#1\rvert}%

\hypersetup{
colorlinks=true,
citecolor=red!70,
linkcolor=green!50!black!70,
urlcolor=red!70
}

\articletype{Research Article}%

\received{26 April 2016}
\revised{6 June 2016}
\accepted{6 June 2016}

\begin{document}

\title{Joint Synthesis of Trajectory and Controlled Invariant Funnel for Discrete-time Systems with Locally Lipschitz Nonlinearities}

\author[1]{Taewan Kim*}

\author[1]{Purnanand Elango}

\author[1]{Beh\c{c}et A\c{c}\i kme\c{s}e}

%

\address[1]{\orgdiv{William E. Boeing Department of Aeronautics and Astronautics}, \orgname{University of Washington}, \orgaddress{\state{Washington}, \country{United States}}}



\corres{*Taewan Kim, \email{twankim@uw.edu}}


\abstract[Summary]{This paper presents a joint synthesis algorithm of trajectory and controlled invariant funnel (CIF) for locally Lipschitz nonlinear systems subject to bounded disturbances. The CIF synthesis refers to a procedure of computing controlled invariance sets and corresponding feedback gains. In contrast to existing CIF synthesis methods that compute the CIF with a pre-defined nominal trajectory, our work aims to optimize the nominal trajectory and the CIF jointly to satisfy feasibility conditions without the relaxation of constraints and obtain a more cost-optimal nominal trajectory. The proposed work has a recursive scheme that mainly optimize trajectory update and funnel update. The trajectory update step optimizes the nominal trajectory while ensuring the feasibility of the CIF. Then, the funnel update step computes the funnel around the nominal trajectory so that the CIF guarantees an invariance property. As a result, with the optimized trajectory and CIF, any resulting trajectory propagated from an initial set by the control law with the computed feedback gain remains within the feasible region around the nominal trajectory under the presence of bounded disturbances. \change{We validate the proposed method via two applications from robotics.}}
 


\maketitle


\section{Introduction}

\change{There has been significant amount of research on trajectory planning algorithms for nonlinear dynamics with nonconvex constraints \cite{malyuta2022convex}. A primary challenge in this area is handling uncertainty such as external disturbances, model mismatch in system dynamics, and state estimation error. Relying solely on a generated nominal trajectory, consisting of time-varying state and open-loop input signals, may not yield a robust control system. This is because the uncertainty can cause the system to deviate from the nominal trajectory. A potential solution is to synthesize not only the open-loop control but also feedback control law, which prevents the system from deviating too far from the nominal trajectory in the presence of uncertainties. To ensure safety by satisfying the prescribed constraints, it is necessary to compute a forward invariant or reachable set, also referred to as a controlled invariant funnel (CIF), that encapsulates all potential state trajectory under the uncertainty. Then, ensuring that the CIF remains inside the safety region can guarantee that all potential trajectories of the system, starting from the funnel entry, will remain safe. Consequently, there has been active research aimed at optimizing both the feedback control and the CIF together in conjunction with the safety constraints. This process of computing the CIF and the associated feedback controller is often referred to as funnel synthesis \cite{reynolds2021funnel,10167750}.}

\begin{figure}
    \centering
        \centering
         \begin{subfigure}[b]{0.3\textwidth}
         \caption{}
         \centering
         \includegraphics[width=0.8\textwidth,trim={2cm 1.5cm 1.5cm 2cm},clip]{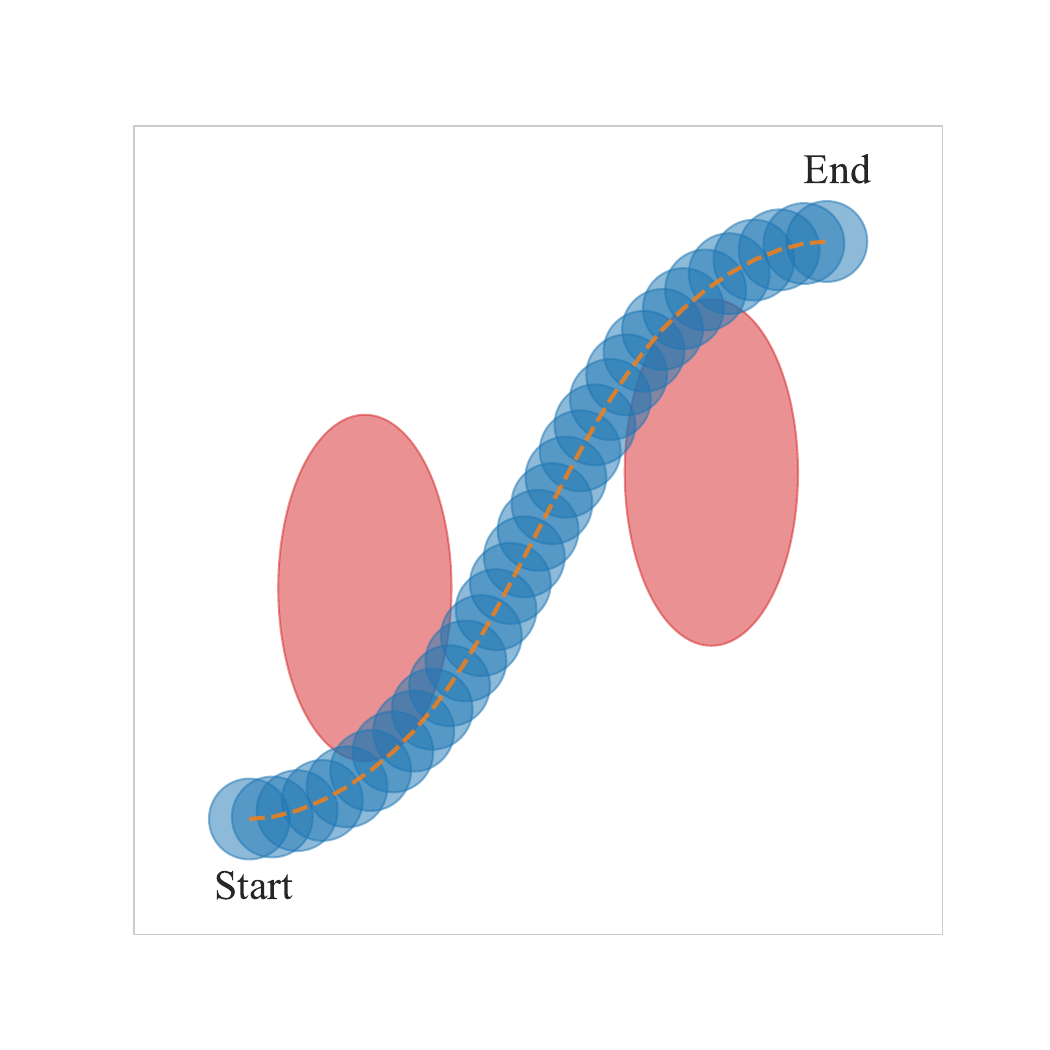}
         \label{fig:motivation-a}
         \end{subfigure}
         \begin{subfigure}[b]{0.3\textwidth}
         \caption{}
         \centering
         \includegraphics[width=0.8\textwidth,trim={2cm 1.5cm 1.5cm 2cm},clip]{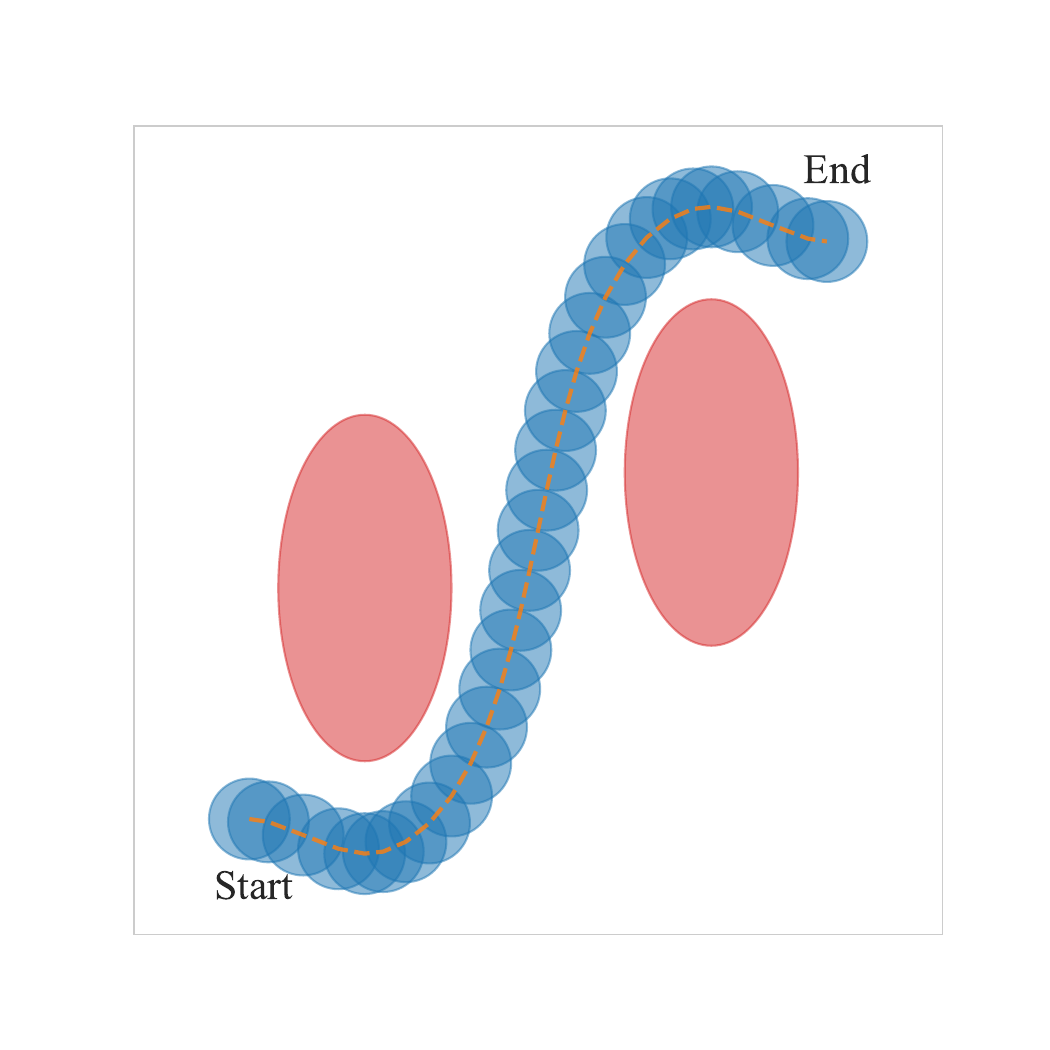}
         \label{fig:motivation-b}
         \end{subfigure}
         \begin{subfigure}[b]{0.3\textwidth}
         \caption{}
         \centering
         \includegraphics[width=0.8\textwidth,trim={2cm 1.5cm 1.5cm 2cm},clip]{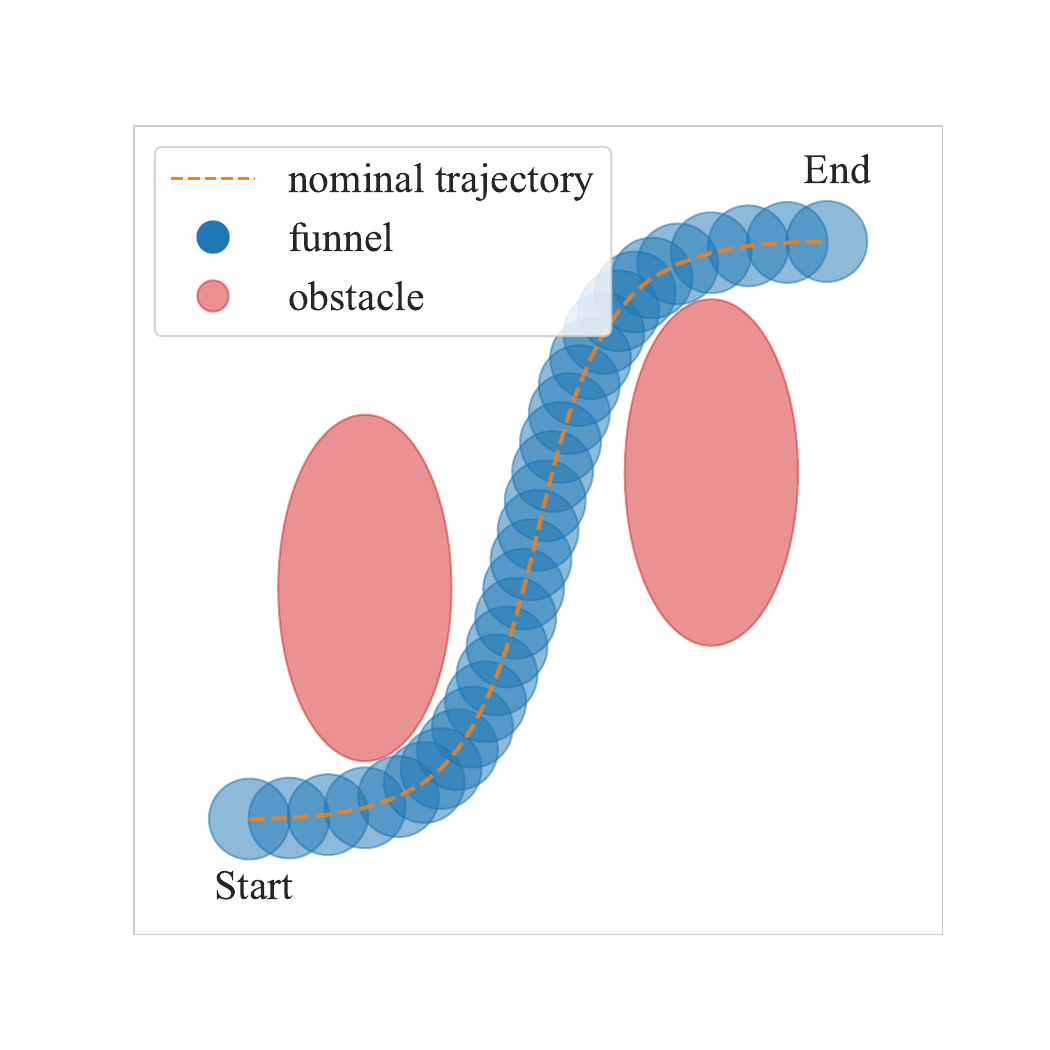}
         \label{fig:motivation-c}
         \end{subfigure}
\caption{\change{Comparative illustration of separate synthesis and joint synthesis approahces. Left: the funnel exhibits constraint violation due to an underestimated safety margin for uncertainties. Middle: the trajectory and the funnel are feasible, but suboptimal due to an overestimated safety margin for uncertainties. Right: the proposed joint synthesis yields an optimal trajectory and funnel while satisfying the obstacle avoidance constraint.}}
\label{fig:motivation}
\end{figure}

\change{The design of robust controllers for uncertain systems traditionally follows a two-step (or separate) scheme; initially, nominal trajectory planning \cite{bonalli2019gusto,malyuta2022convex} is performed to compute the open-loop control and the corresponding state trajectory, followed by the synthesis of the feedback control and the associated CIF \cite{majumdar2017funnel,buch2021finite} based on the analysis of the perturbed system around the nominal trajectory. In aerospace applications, the former is often termed to as \emph{guidance} and the latter is referred to as \emph{control}. However, such a two step scheme has a potential drawback; the resulting control law, consisting of both the open-loop and closed-loop control, may be overly conservative. This conservatism stems from the lack of joint consideration of the open-loop and closed-loop control computations for given constraints such as actuator limits and obstacles. Consider, for example, a path planning scenario shown in Fig.~\ref{fig:motivation} where there are two obstacles between the start and end points. If the nominal trajectory is optimized independently of the CIF (Fig. \ref{fig:motivation-a}), the resulting trajectory may be close to the obstacle boundary, so that the trajectory cost such as minimum-time or minimum-fuel is optimized. In such cases, the CIF can violate the constraints because the CIF size cannot be arbitrarily minimized under the presence of the external disturbances. One way to resolve issue is to introduce a safety margin in the constraints for the nominal trajectory planning. However, since it is not tractable to know the optimal margin beforehand, if the margin is too large, the resulting nominal trajectory and CIF can be conservative (Fig. \ref{fig:motivation-b}).}

\change{To overcome this drawback, we propose an algorithm referred to as a \emph{joint synthesis}, that performs the funnel synthesis jointly with the nominal trajectory optimization in an iterative scheme. The key advantage of the proposed approach for the joint synthesis of trajectory and CIF is its ability to reduce conservatism that can be caused with the aforementioned separate synthesis, and hence improve the optimality of the resulting control design without compromising robustness. The proposed joint synthesis algorithm does not need to estimate the safety margin of the given uncertainty in advance. Hence, this allows the algorithm to exploit a larger feasible set while designing the trajectory and the CIF jointly, thereby generating a more optimal trajectory and CIF (Fig.~\ref{fig:motivation-a}).}

Specifically, the proposed method jointly generates trajectory and synthesizes funnel for locally Lipschitz nonlinear systems subject to bounded disturbances. The problem formulation can be viewed as a robust trajectory optimization in which we optimize both the trajectory and the CIF that consists of the forward invariant set and the corresponding feedback gain. To this end, we draw ideas from sequential convex programming (SCP) \cite{bonalli2019gusto}, Lyapunov theory, and linear matrix inequalities (LMIs) for robust control \cite{xu2020observer}. The proposed method has the following steps in each iteration: First, we update the nominal trajectory while ensuring the feasibility of the funnel. The next step estimates local Lipschitz constants of the nonlinearity in the system by sampling state space inside the funnel. With the trajectory computed in the first step, the third step then constructs a semidefinite programming (SDP) problem derived with \change{funnel constraints} and a Lyapunov condition to ensure the invariance property of the funnel. These steps are repeated until the convergence of both the trajectory and the funnel synthesis. We validate the proposed method through a numerical simulation.

\subsection{Related work}
\change{In this subsection, we discuss related work with the comparison with the proposed work.} 
Optimizing the trajectory and CIF jointly has been studied in the context of robust MPC \cite{bemporad1999robust}. To reduce the computational complexity to satisfy the real-time performance, most works in robust MPC precompute the feedback gain or the forward reachable set and then optimize the nominal trajectory online. For example, a tube-based MPC scheme is developed for Lipschitz nonlinear systems subject to bounded disturbances in \citenum{yu2013tube}. With the precomputed feedback gain and the tube set, the proposed method optimizes the nominal trajectory. The work in \citenum{kohler2020computationally} computes an incremental Lyapunov function and a corresponding feedback gain offline, and then optimizes the nominal trajectory and the support value of the invariant set online. 
In \citenum{villanueva2017robust}, they obtain both the nominal trajectory and the invariant set by solving an SDP with the precomputed feedback gain. On the other hand, the proposed work separately parameterizes the nominal trajectory and the CIF that consists of the invariant state set and the feedback controller, and then optimizes them together in the recursive scheme. For nonlinear systems having incrementally conic uncertainties/nonlinearities, \citenum{accikmecse2011robust} provides an LMI-based framework for the generation of the control policy and the invariant set for robust MPC. Our proposed work can be viewed as an extension of \citenum{accikmecse2011robust} in order to handle more general nonlinear systems that are locally Lipschitz.

\change{The work in \citenum{garimella2018robust} proposes a robust tube-based MPC algorithm for obstacle avoidance. They provide a framework that optimizes the nominal trajectory and the corresponding funnel based on the derived ellipsoid uncertainty propagation. However, this propagation provides only an approximate guarantee, that is, the resulting funnel might not be invariant under the bounded disturbance, and they assume that the feedback controllers are a priori given. The proposed method, however, provides the exact invariance guarantee and optimizes the controller simultaneously, thereby reducing conservatism. For further examples of obstacle avoidance using MPC, we refer the reader to \citenum{richards2006robust,di2012model,seo2022real}.} The most relevant results in robust MPC literature to the proposed work appear in \citenum{messerer2021efficient} and \citenum{manchester2019robust}, which jointly synthesize the nominal trajectory and the funnel. The major difference compared to our method is that the resulting funnel in \citenum{messerer2021efficient} and \citenum{manchester2019robust} is only robust for the linearized closed-loop system (they ignored the higher-order terms due to nonlinearities). On the other hand, the proposed method considers the higher-order terms as state and input-dependent uncertainty using the local Lipschitz property. Thus, the robustness of the resulting funnel is not compromised. Another difference is that the proposed method sets the linear feedback gains and the invariant set variables as decision variables and then optimizes them simultaneously to satisfy the invariance and feasibility properties of the funnel. In contrast, the invariant set variables are determined by the ellipse propagation equation derived from the linearized closed-loop system in \citenum{messerer2021efficient}. In \citenum{manchester2019robust}, the feedback gain is computed by solving the time-varying linear quadratic regulator (LQR) problem with the linearized closed-loop system. Again, in order to do that, both utilized the linearized dynamics.

\change{The CIF computation in the proposed work is relevant to recent studies on finite horizon robustness analysis \cite{seiler2019finite} and robust synthesis \cite{buch2021finite}. The finite horizon robust analysis for uncertain linear time-varying (LTV) systems is studied in \citenum{seiler2019finite}, where the behavior of the uncertain system is characterized by integral quadratic constraints (IQCs). This work is extended in \citenum{buch2021finite} where the robust controller is synthesized based on the established robustness analysis. Similar to \citenum{seiler2019finite} and \citenum{buch2021finite}, we construct the CIF of the nonlinear systems by describing the nonlinearity with (incremental) quadratic inequality that can be viewed as a pointwise (incremental) IQC. The major difference is that the approach in \citenum{seiler2019finite,buch2021finite} focuses on analyzing and synthesizing the LTV systems rather nonlinear systems. To apply their approaches for the nonlinear systems, one needs obtain such LTV systems via linearization around a given nominal trajectory, but obtaining such nominal trajectories is not considered in their work. Meanwhile, the proposed work indeed considers this; we design an algorithm that optimizes the nominal trajectory in conjunction with synthesizing the uncertain LTV system obtained via linearization, computing both open-loop and feedback controllers, and the reachable set (funnel) for the nonlinear systems.}

The proposed work is motivated by studies on the robust CIF generation. In \citenum{majumdar2017funnel}, sum-of-squares (SOS) programming is applied to design the CIF for nonlinear systems having polynomial dynamics subject to disturbances. They first design a finite library of open-loop nominal trajectories. Then, for each nominal trajectory, they optimize a feedback controller and a corresponding invariant set by solving SOS programming iteratively. Hence, their method is one shot procedure where the nominal trajectory is computed first, and then the computation of the funnel follows. The similar one shot approaches are conducted in \citenum{seo2021fast,jang2021fast}. For the fast computation of the CIF, the work in \citenum{seo2021fast} formulates an optimization problem for establishing the CIF as a linear program (LP) which is computationally cheaper than SOS programming. This research is extended to consider piecewise polynomial systems in \citenum{jang2021fast}. These work, however, do not consider the controller synthesis, and hence focus on obtaining the reachable set (funnel) of the given (polynomial) closed-loop system. \change{An LMI-based CIF computation for locally Lipschitz nonlinear systems is studied in our previous work \citenum{reynolds2021funnel} where the CIF computation is formulated as SDP that is also computationally efficient to solve than the problem with SOS programming, and the work is extended in \citenum{10167750} to consider the bounded external disturbance. The aforementioned work including our previous work solely focuses on generating the CIF given the nominal trajectory, whereas the proposed work has a recursive scheme where the nominal trajectory, the feedback controller, and the invariant sets are optimized simultaenously.}

\change{Recently, a robust control method that jointly optimizes the nominal trajectory and the linear feedback gain was studied in \citenum{leeman2023robust1}. Using a system level synthesis framework, they formulated the problem of obtaining the nominal trajectory and the linear feedback gain for the LTV systems describing the error dynamics around the nominal trajectory. This work is extended to handle parametric uncertainties in \citenum{leeman2023robust2}. The advantage of the proposed work is that it does not assume the disturbance to be additive, whereas in \citenum{leeman2023robust1,leeman2023robust2} the bounded disturbances simply is added to the system dynamics. Hence, the proposed method can handle cases where the bounded disturbances manifest via a nonlinear relationship. We demonstrate this capability with the numerical simulation. Additionally, the approach in \citenum{leeman2023robust1,leeman2023robust2} assumes system dynamics to be three times continuously differentiable, whereas the proposed work assumes one-time continuous differentiability.}

\subsection{Contributions}

\change{The contributions of the proposed work can be summarized as follows: First, we propose a novel algorithm that jointly synthesizes the nominal trajectory and the CIF to ensure robustness for nonlinear systems with locally Lipschitz nonlinearities. By jointly optimizing them together, the algorithm can mitigate the potential conservatism that may arise from optimizing them separately. Second, the proposed method can be applicable to a broad class of nonlinear systems that have locally Lipschitz nonlinearities (including systems with continuously differentiable dynamics) and non-additive disturbances.  Third, we extend the existing LMI-based CIF computation and robust MPC research \cite{accikmecse2011robust,reynolds2021funnel,10167750} in order to handle more general class of nonlinear systems that are discrete-time locally Lipschitz systems under the presence of norm-bounded uncertainties. Finally, we validate the proposed method through two distinct robotic applications: the first involves path planning with obstacle avoidance for a unicycle model, and the second focuses on a 6-degree-of-freedom (6-DoF) free-flying spacecraft.}

\subsection{Outline}
We present the problem formulation in Section \ref{sec2} and the proposed method in Section \ref{sec3}. In Section \ref{sec4}, we perform a numerical evaluation for our proposed method. Concluding remarks are provided in \ref{sec5}.

\subsection{Notation}

Let $\mathbb{R}$ be the field of real numbers, $\mathbb{R}^{n}$ be the $n$-dimensional Euclidean space, and $\mathbb{N}$ be the set of natural numbers. A finite set of consecutive non-negative integers is represented by $\mathcal{N}_{q}^{r}\coloneqq\{q,q+1,\ldots,r\}$. The symmetric matrix $Q=Q^{\top}(\succeq)\!\succ 0$ implies $Q$ is positive-(semi-)definite matrix, and $(\mathbb{S}^n_{+})\, \mathbb{S}^n_{++}$ denotes the set of all positive-(semi-)definite matrices whose size is $n\times n$. The symbol $\oplus$ denotes the Minkowski sum. The vector $(x,y)$ represents the concatenation of two vectors $x$ and $u$ into a longer vector. The notation {*} represents the symmetric part of a matrix, i.e, $\left[\begin{array}{cc}
a & b^{\top}\\
b & c
\end{array}\right]=\left[\begin{array}{cc}
a & *\\
b & c
\end{array}\right]$, and $\{\bar{x}_{k},\bar{u}_{k},\bar{w}_{k}\}_{k=0}^{K}$ illustrates  $\{\bar{x}_{0},\bar{u}_{0},\bar{w}_{0},\ldots,\bar{x}_{K},\bar{u}_{K},\bar{w}_{K}\}$. The symmetric squared root of a symmetric matrix $A$ is defined as $A^\frac{1}{2}$ by eigenvalue decomposition \cite{boyd2004convex}. The operation $\text{diag}(\cdot)$ is a diagonal matrix formed from its vector argument.

\section{Problem formulation}
\label{sec2}
Consider a discrete-time uncertain nonlinear system of the following
form:
\begin{align}
x_{k+1}=f(t_{k},x_{k},u_{k},w_{k}),\quad\forall\, k\in\mathcal{N}_{0}^{N-1},\label{eq:nonlinear}
\end{align}
where $N\in\mathbb{N}$ is the length of the time horizon and $t_k\in\mathbb{R}$ is the time at the $k$. The function $f:\mathbb{R}\times\mathbb{R}^{n_{x}}\times\mathbb{R}^{n_{u}}\times\mathbb{R}^{n_{w}}\rightarrow\mathbb{R}^{n_{x}}$ is assumed to be a locally Lipschitz and at least once differentiable. The vector $x_{k}\in\mathbb{R}^{n_{x}}$ is the state, $u_{k}\in\mathbb{R}^{n_{u}}$ is the control input, and the signal $w_{k}\in\mathbb{R}^{n_{w}}$ is the exogenous disturbance or model mismatch that is assumed to be unknown but norm bounded: $\norm{w_{k}}_{2}\leq1$ for all $k\in\mathcal{N}_{0}^{N-1}$. 

Let $\{\bar{x}_{k}\}_{k=0}^N, \{\bar{u}_{k},\bar{w}_k\}_{k=0}^{N-1}$ be a nominal trajectory that the CIF is centered around, and is feasible for the nonlinear dynamics \eqref{eq:nonlinear}. In this paper, the nominal trajectory is assumed to have zero disturbances, i.e., $\bar{w}_{k}=0$ for all $k\in\mathcal{N}_{0}^{N-1}$. We define difference state $\eta_{k}\coloneqq x_{k}-\bar{x}_{k}$ and difference input $\xi_{k}\coloneqq u_{k}-\bar{u}_{k}$, and assume a linear feedback $\xi_{k}=K_{k}\eta_{k}$ for all $k\in\mathcal{N}_{0}^{N-1}$, which leads to a closed-loop system and a control law given by
\begin{align}
\eta_{k+1} & =f(t_{k},x_{k},u_{k},w_{k})-f(t_{k},\bar{x}_{k},\bar{u}_{k},0),\label{eq:closed_loop_general}\\
u_{k}&=\bar{u}_{k}+K_{k}\eta_k,\quad\forall\, k\in\mathcal{N}_{0}^{N-1},   \label{eq:control_law}
\end{align}
where $K_{k}\in\mathbb{R}^{n_{x}\times n_{u}}$ is a feedback gain. In this paper, we consider a specific class of funnels that consists of ellipsoids of state and input. The ellipsoid for the difference state is represented as
\begin{align}
\mathcal{E}_{Q_{k}}\coloneqq\{\eta\in\mathbb{R}^{n_{x}}\mid\eta^{\top}Q_{k}^{-1}\eta\leq1\},\quad\forall\, k\in\mathcal{N}_{0}^{N},\label{eq:def_ellipse}
\end{align}
where $Q_{k}\in\mathbb{S}_{++}^{n_{x}\times n_{x}}$ is a positive definite matrix. With the linear feedback gain $K_{k}$, it follows from Schur complement that $\eta_{k}\in\mathcal{E}_{Q_{k}}$ implies $\xi_{k}\in\mathcal{E}_{K_{k}Q_{k}K_{k}^\top}$ \cite{reynolds2021funnel}. Now we are ready to formally define the quadratic CIF.
\begin{definition}
A quadratic controlled positively invariant funnel, $\mathcal{F}_{k}$, associated with a closed loop system \eqref{eq:closed_loop_general} is a time-varying set in state and control space that is parameterized by a time-varying positive definite matrix $Q_{k}\in\mathbb{S}_{++}^{n_{x}}$ and a time-varying matrix $K_{k}\in\mathbb{R}^{n_{x}\times n_{u}}$ such that $\mathcal{F}_{k}=\mathcal{E}_{Q_{k}}\times\mathcal{E}_{K_{k}Q_{k}K_{k}^{\top}}$, and the funnel $\mathcal{F}_{k}$ is invariant and 
feasible for all $k\in \mathcal{N}_0^{N}$.
\end{definition}
The invariance property of the CIF with the closed-loop system \eqref{eq:closed_loop_general} and the control law \eqref{eq:control_law} can be mathematically stated as follows:
\begin{align}
(\eta_{0},\xi_{0}) & \in\mathcal{F}_{0}\Rightarrow(\eta_{k},\xi_{k})\in\mathcal{F}_{k},\:\forall\, k\in\mathcal{N}_{1}^{N}.\label{eq:funnel_invariant}
\end{align}
This condition implies that if a particular initial condition is inside the funnel, then a trajectory propagated with the closed-loop model
\eqref{eq:closed_loop_general} remains within the funnel as well. The feasibility property for the funnel $\mathcal{F}_{k}$ can be mathematically expressed as:
\begin{subequations}
\label{eq:funnel_feasibility}
\begin{align}
\{\bar{x}_{k}\}\oplus\mathcal{E}_{Q_{k}} & \subseteq\mathcal{X},\\
\{\bar{u}_{k}\}\oplus\mathcal{E}_{K_{k}Q_{k}K_{k}^\top} & \subseteq\mathcal{U},\quad\forall\, k\in\mathcal{N}_{0}^{N-1}.
\end{align}
\end{subequations}
The feasibility conditions require that every state and input in the funnel around the nominal trajectory should be feasible for the given state and input constraint sets $\mathcal{X}$ and $\mathcal{U}$, respectively.

Now we are ready to derive the problem formulation. The goal of the joint synthesis of trajectory and CIF is to solve a discrete-time nonconvex optimization problem of the following form:
\begin{subequations}\label{eq:main_all}
\begin{align}
\underset{
\begin{array}{c} \scriptscriptstyle \bar{x}_{k},Q_{k},{\mu^Q_{k}}\,\forall\, k\in\mathcal{N}_{0}^{N}, \\[-0.1cm]
\scriptscriptstyle \bar{u}_{k},K_{k},{\mu^K_{k}}, \forall\, k\in\mathcal{N}_{0}^{N-1}
\end{array}}{\operatorname{minimize}} & \change{\sum_{k=0}^{N-1} J_{t}(\bar{x}_{k},\bar{u}_{k})
+ w_Q\sum_{k=0}^{N} \mu^Q_{k} + w_K\sum_{k=0}^{N-1} \mu^K_{k}} \label{eq:maincost}\\[-0.1cm]
\operatorname{subject~to}~~&
\bar{x}_{k+1}=f(t_{k},\bar{x}_{k},\bar{u}_{k},0),\forall\, k\in\mathcal{N}_{0}^{N-1} \label{eq:main_dynamics}\\
 & Q_{k}\preceq{\mu^Q_{k}}I,\forall\, k\in\mathcal{N}_{0}^{N}\label{eq:main_minimizeQ}\\
 & K_{k}Q_{k}K_{k}^{\top}\preceq{\mu^K_{k}}I,\forall\, k\in\mathcal{N}_{0}^{N-1}\label{eq:main_minimizeK}\\
 & \operatorname{conditions~} \eqref{eq:funnel_invariant}-\eqref{eq:funnel_feasibility},\nonumber \\
 & \bar{x}_{0}\oplus\mathcal{E}_{Q_{0}}\supseteq\mathcal{X}_{0},\\
 & \bar{x}_{N}\oplus\mathcal{E}_{Q_{N}}\subseteq\mathcal{X}_{N},\label{eq:mainboundary}
\end{align}
\end{subequations}
where the summands in the objective function consist of the trajectory cost and the funnel cost, \change{and $0<w_Q\in\mathbb{R}$ and $0<w_K\in\mathbb{R}$ are user-defined weights.} The function $J_{t}$ is a cost for the trajectory and is assumed to be convex in $\bar{x}_{k}$ and $\bar{u}_{k}$. The slack variables ${\mu^Q_{k}}\in\mathbb{R}$ and ${\mu^K_{k}}\in\mathbb{R}$ are introduced to minimize the diameter of the ellipsoidal sets $\mathcal{E}_{Q_{k}}$ and $\mathcal{E}_{K_{k}Q_{k}K_{k}^{\top}}$ in the funnel by imposing the constraints \eqref{eq:main_minimizeQ}-\eqref{eq:main_minimizeK}. Minimizing the size of the funnel leads the effect of the propagated disturbances starting from the initial set to be minimized \cite{majumdar2017funnel}. 
While minimizing the cost, the formulation guarantees the invariance property in \eqref{eq:funnel_invariant} and ensures the feasibility of the ellipsoids encapsulating the nominal states and inputs in \eqref{eq:funnel_feasibility}. For boundary conditions, the initial and final ellipsoids, $\mathcal{X}_{0}$ and $\mathcal{X}_{N}$, are given as
\begin{subequations}
\label{eq:boundary_details}
\begin{align}
\mathcal{X}_{0} & =\big\{x\,\big|\,(x-x_{i})^{\top}Q_{i}^{-1}(x-x_{i})\leq1\big\},\\
\mathcal{X}_{N} & =\big\{x\,\big|\,(x-x_{f})^{\top}Q_{f}^{-1}(x-x_{f})\leq1\big\}, 
\end{align}
\end{subequations}
where $x_{i}\in\mathbb{R}^{n_{x}}$ is a nominal initial state, $Q_{i}\in\mathbb{S}_{++}^{n_{x}}$ is a constant matrix defining the initial ellipsoidal set, $x_{f}\in\mathbb{R}^{n_{x}}$ is the nominal final state, and $Q_{f}\in\mathbb{S}_{++}^{n_{x}}$ is a constant matrix defining the final ellipsoidal set. The computed funnel at $k=0$ should include the initial set $\mathcal{X}_{0}$ to generate the trajectory from any state in the initial set. Also, the ellipsoid corresponding to the state in the funnel at $k=N$ should be a subset of $\mathcal{X}_{N}$ so that the resulting trajectory is guaranteed to terminate in $\mathcal{X}_N$.

\begin{figure}
\begin{center}
\includegraphics[width=8cm]{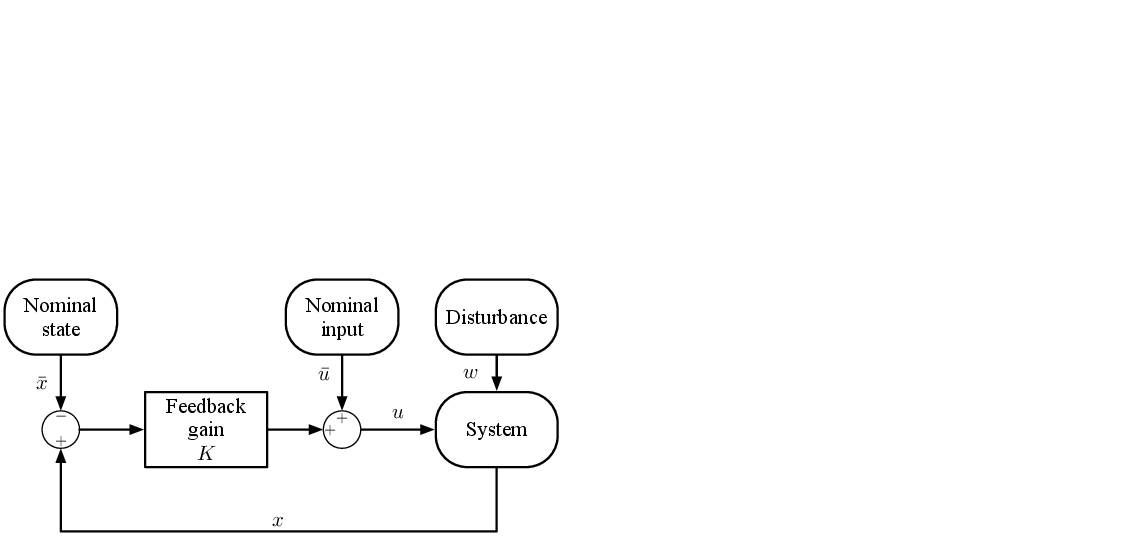}  
\caption{
\change{A block diagram of the control procedure.}
} 
\label{fig:controldiagram}
\end{center}
\end{figure}

It is worth mentioning that the system dynamics \eqref{eq:main_dynamics} for the nominal trajectory has no disturbances ($\bar{w}_{k}=0$), but the invariance property is achieved with the closed-loop dynamics \eqref{eq:closed_loop_general}-\eqref{eq:control_law} in which the disturbances exist. Hence, any trajectory propagated for the uncertain nonlinear dynamics \eqref{eq:nonlinear} with the control law \eqref{eq:control_law} from any initial state in $\mathcal{X}_{0}$ remains within the feasible region under the presence of norm bounded uncertainties. \change{The block diagram of the resulting control signal is illustrated in Fig.~\ref{fig:controldiagram}.}

\section{Iterative Robust Trajectory Optimization}
\label{sec3}
\begin{figure}
\begin{center}
\includegraphics[width=9cm]{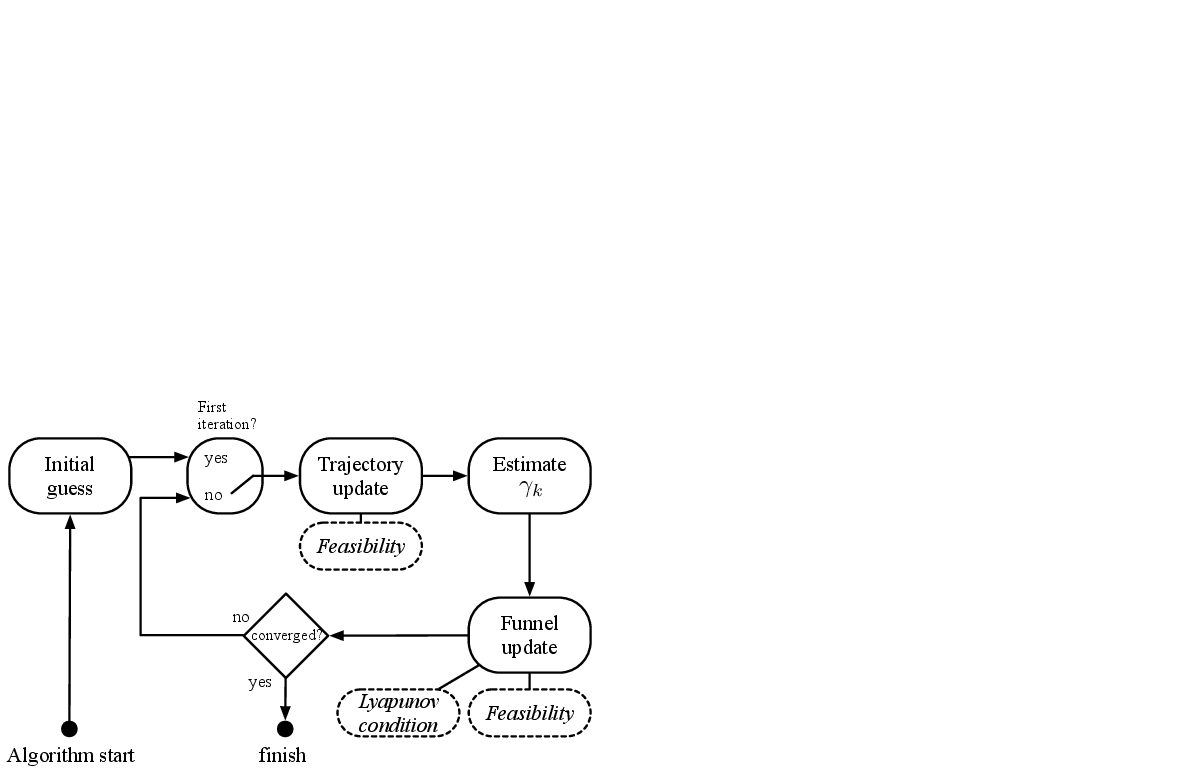}  
\caption{
A block diagram of the proposed method. Starting from the initial guess, the method optimizes the trajectory while considering the feasibility of the funnel. The local Lipschitz constant $\gamma_k$ of the nonlinearity around the obtained trajectory is then estimated. The next step is to optimize the funnel with \change{the funnel constraints and} the Lyapunov condition that ensures the invariance property. The entire process is repeated until both the trajectory and the funnel converge.
} 
\label{fig:blockdiagram}
\end{center}
\end{figure}

In this section, we discuss the details of the proposed method to solve the robust trajectory optimization problem given in \eqref{eq:maincost}-\eqref{eq:mainboundary}. The method tackles the problem by iteratively updating the nominal trajectory $\{\bar{x}_k\}_{k=0}^{N}$,$\{\bar{u}_k\}_{k=0}^{N-1}$, the parameters of the set $\{Q_{k}\}_{k=0}^{N}$ and the feedback controller $\{K_{k}\}_{k=0}^{N}$ in the CIF. In each iteration, the method consists of 3 steps: the nominal trajectory update, the estimation of the locally Lipschitz constant, and the funnel update. In this section, we denote an initial guess or solution variables of the previous iteration (i.e., reference trajectory and funnel parameters) by $\{\hat{x}_{k},\hat{Q}_{k}\}_{k=0}^{N}, \{\hat{u}_{k},\hat{K}_{k}\}_{k=0}^{N-1}$. The block diagram of the proposed algorithms is given in Fig.~\ref{fig:blockdiagram}. 

\subsection{Nominal trajectory update}
We require the nominal trajectory to satisfy the (possibly nonconvex) constraints \eqref{eq:main_dynamics} and \eqref{eq:funnel_feasibility} while minimizing the trajectory cost $J_{t}$ by approximating the original problem with a convex sub-problem. This is a typical process in many SCP methods to solve nonconvex trajectory optimization problems \cite{malyuta2022convex}. In contrast to the typical SCP methods, the feasibility problem in \eqref{eq:funnel_feasibility} involves the funnel parameters that are fixed as the reference funnel variables $\{\hat{Q}_{k},\hat{K}_{k}\}_{k=0}^{N-1}$ in this trajectory update step.

In each sub-problem, the intermediate trajectory solution should
satisfy the following affine system:
\begin{align}
x_{k+1}=\bar{A}_{k}x_{k}+\bar{B}_{k}u_{k}+\bar{z}_{k}+v_{k},\quad\forall\, k\in\mathcal{N}_{0}^{N-1}\label{eq:lin_dynamics}
\end{align}
where $\bar{A}_{k},\bar{B}_{k},\bar{z}_{k}$ define the linearized model of the nonlinear dynamics given in \eqref{eq:nonlinear} evaluated around the reference trajectory $\{\hat{x}_{k},\hat{u}_{k}\}_{k=0}^{N-1}$ with zero disturbance $\bar{w}_{k}=0$. The term $v_{k}$ is a virtual control variable that serves to prevent the sub-problem from being artificially infeasible \cite{malyuta2022convex} due to linearization of dynamics and constraints.

The feasible sets $\mathcal{X}$ and $\mathcal{U}$ are expressed as
\begin{align*}
\mathcal{X} & =\{x\mid h_{i}(x)\leq0,\quad i=1,\ldots,m_{x}\},\\
\mathcal{U} & =\{u\mid {g_j}(u)\leq0,\quad j=1,\ldots,m_{u}\},
\end{align*}
where $h_{i}$ and ${g_j}$ are at least once differentiable functions. \change{While we assume here that $\mathcal{X}$ and $\mathcal{U}$ are time-invariant for brevity, the proposed framework, however, can easily incorporate time-varying sets.} The nonlinear constraints need to be linearized to ensure convexity of the sub-problem. Thus, we approximate the feasible set $\mathcal{X}$ and $\mathcal{U}$ as polytopes, which are obtained via linearization around $\{\hat{x}_{k},\hat{u}_{k}\}_{k=0}^{N-1}$ as follows:
\begin{align*}
\mathcal{P}_{k}^{x} & =\{x\mid(a_{i}^{x})_{k}^{\top}x_{k}\leq(b_{i}^{x})_{k},\quad i=1,\ldots,m_{x}\},\\
\mathcal{P}_{k}^{u} & =\{u\mid(a_{j}^{u})_{k}^{\top}x_{k}\leq(b_{j}^{u})_{k},\quad j=1,\ldots,m_{u}\},
\end{align*}
where $(a_{i}^{x},b_{i}^{x})$ and $(a_{j}^{u},b_{j}^{u})$ are first-order approximations of $h_{i}$ and ${g_j}$, respectively. \change{Notice that while $\mathcal{X}$ and $\mathcal{U}$ are assumed to be time-invariant, their polytopic approximations $\mathcal{P}_k^x$ and  $\mathcal{P}_k^u$ would be time-varying due to the time variation in the reference trajectory $\hat{x}_k, \hat{y}_k$.} Then, the feasibility conditions with the fixed funnel parameters $\{\hat{Q}_{k},\hat{K}_{k}\}_{k=0}^{N-1}$ in \eqref{eq:funnel_feasibility} can be approximated as linear constraints as follows \cite{reynolds2021funnel}:
\begin{subequations}
\label{eq:lin_feasibility}
\begin{align}
\norm{(\hat{Q}_{k}^{\top})^{\frac{1}{2}}(a_{i}^{x})_{k}}_{2}+(a_{i}^{x})_{k}^{\top}x_{k} & \leq(b_{i}^{x})_{k},i=1,\ldots,m_{x}, \\
\norm{(\hat{K}_{k}\hat{Q}_{k}\hat{K}_{k}^{\top})^{\frac{1}{2}}(a_{j}^{u})_{k}}_{2}+(a_{j}^{u})_{k}^{\top}u_{k} & \leq(b_{j}^{u})_{k},j=1,\ldots,m_{u}, \quad \forall\, k\in\mathcal{N}_{0}^{N-1}.
\end{align}
\end{subequations}

The trajectory update step for the nominal trajectory has the following
form of a second-order cone program (SOCP):
\begin{subequations}\label{eq:traj_update}
\begin{align}
\underset{\begin{array}{c}
\scriptscriptstyle \bar{x}_{k},\bar{u}_{k},\bar{v}_{k},\bar{x}_{N},\\[-0.1cm]
\scriptscriptstyle \forall\, k\in\mathcal{N}_{0}^{N-1}
\end{array}}{\operatorname{minimize}} & \sum_{k=0}^{N-1}J_{t}(\bar{x}_{k},\bar{u}_{k})+J_{vc}(\bar{v}_{k})+J_{tr}(\bar{x}_{k},\bar{u}_{k})\label{eq:n_cost}\\[-0.25cm]
\operatorname{subject~to}~& \operatorname{conditions~} \eqref{eq:lin_dynamics}-\eqref{eq:lin_feasibility},\nonumber \\
 & \bar{x}_{0}=x_{i},\quad \bar{x}_{N}=x_{f}.\label{eq:n_boundary}
\end{align}
\end{subequations}
In the cost function, there are two additional penalty terms for virtual control $J_{vc}$ and trust region $J_{tr}$. The virtual control penalty
enforces the virtual control variables $v_{k}$ to remain small, and the trust region encourages the optimum to stay in the vicinity of the reference trajectory $\{\hat{x}_{k},\hat{u}_{k}\}_{k=0}^{N-1}$ where the linearization error is small. They are formulated as follows:
\begin{subequations}
\label{eq:trust_region_virtual_control}
\begin{align}
J_{vc}(v_{k}) & =w_{v}\norm{v_{k}}_{1}, \\
J_{tr}(x_{k},u_{k}) & =w_{tr}(\norm{x_{k}-\hat{x}_{k}}_{2}^{2}+\norm{u_{k}-\hat{u}_{k}}_{2}^{2}), 
\end{align}
\end{subequations}
where $w_{v}\in\mathbb{R}$ and $w_{tr}\in\mathbb{R}$ are user-defined weight parameters for the virtual control and the trust region, respectively. As a result of the optimization problem \eqref{eq:traj_update}, the solution becomes a new nominal trajectory $\{\bar{x}_{k}\}_{k=0}^{N},\{\bar{u}_{k}\}_{k=0}^{N-1} $ that will be used for the funnel computation in the following section. This type of penalized-trust region-based optimization has been studied for trajectory optimization \cite{reynolds2020dual} and general nonlinear programming \cite{cartis2011evaluation}.

\subsection{CIF update}

In this section, we describe how to optimize the CIF around the nominal trajectory obtained from the previous section. The optimization problem derived in this section aims to make the
funnel invariant \eqref{eq:funnel_invariant} and feasible for \change{the constraints \eqref{eq:funnel_feasibility} and} the boundary conditions \eqref{eq:mainboundary} for locally Lipschitz nonlinear systems. To this end, we construct a SDP whose solution provides the parameters of the invariant set and the feedback gains $\{Q_{k}\}_{k=0}^{N}, \{K_{k}\}_{k=0}^{N-1}$.

\subsubsection{Nonlinear dynamics}

Since the nonlinear dynamics in \eqref{eq:nonlinear} is at least once differentiable, it can be re-written as 
\begin{subequations}
\label{eq:converted_model}
\begin{align}
x_{k+1} & =f(t_k,x_{k},u_{k},w_{k}), \\
 & =A_{k}x_{k}+B_{k}u_{k}+F_{k}w_{k}+Ep_{k}, \\
p_{k} & ={\phi_{k}}(q_{k}), \\
q_{k} & =Cx_{k}+Du_{k}+Gw_{k}.
\end{align}
\end{subequations}
Notice that all nonlinearities are lumped into a vector $p_{k}\in\mathbb{R}^{n_{p}}$ represented by a function ${\phi_{k}}:\mathbb{R}^{n_{q}}\rightarrow\mathbb{R}^{n_{p}}$ with its argument $q_{k}\in\mathbb{R}^{n_{q}}$. The matrix $E\in \mathbb{R}^{n_x\times n_p}$ is introduced since not all states are affected by the nonlinearities. The matrices $A_{k},B_{k}$ and $F_{k}$ can be arbitrary, but we specifically choose $A_{k}$, $B_{k}$, and $F_{k}$ to be the first order approximation of the nonlinear dynamics $f$ around the nominal trajectory as follows:
\begin{align*}
A_{k} & \coloneqq\frac{\partial f}{\partial x}\bigg|_{x=\bar{x}_{k},u=\bar{u}_{k},w=0},\quad B_{k}\coloneqq\frac{\partial f}{\partial u}\bigg|_{x=\bar{x}_{k},u=\bar{u}_{k},w=0},\\
F_{k} & \coloneqq\frac{\partial f}{\partial w}\bigg|_{x=\bar{x}_{k},u=\bar{u}_{k},w=0},\quad\forall\, k\in\mathcal{N}_{0}^{N-1}.
\end{align*}

With difference state $\eta_{k}$ and input $\xi_{k}$, the difference
dynamics can be derived as
\begin{align*}
x_{k+1}-\bar{x}_{k+1} & =A_{k}\eta_{k}+B_{k}\xi_{k}+F_{k}w_{k}+E(p_{k}-\bar{p}_{k}) +f(t_k,\bar{x}_{k},\bar{u}_{k},0)-\bar{x}_{k+1},
\end{align*}
where $\bar{p}_k = \phi_k (\bar{q}_k)$ and $\bar{q}_k = C\bar{x}_k+D\bar{u}_k$. The term $f(t_k,\bar{x}_{k},\bar{u}_{k},0)-\bar{x}_{k+1}$ on the right hand side exists because of the dynamical error in the intermediate nominal trajectory $\{\bar{x}_{k}\}_{k=0}^{N},\{\bar{u}_{k}\}_{k=0}^{N-1}$. This error is gradually reduced as the iteration proceeds because the trajectory update \eqref{eq:traj_update} ensures that the nominal trajectory becomes dynamically feasible for the entire interval. Thus, we intentionally do not consider this error in the funnel update step since it is sufficient for the funnel to satisfy the invariance and feasibility properties with the converged nominal trajectory that is dynamically feasible. The difference dynamics we consider for the funnel update is consequently written as
\begin{align*}
\eta_{k+1} & =A_{k}\eta_{k}+B_{k}\xi_{k}+F_{k}w_{k}+E\delta p_{k}, \\
\delta p_{k} & ={\phi_{k}}(q_{k})-{\phi_{k}}(\bar{q}_{k}),\\
\delta q_{k} & =C\eta_{k}+D\xi_{k}+Gw_{k},
\end{align*}
where $\delta p_{k}\coloneqq p_{k}-\bar{p}_{k}$ and $\delta q_{k}\coloneqq q_{k}-\bar{q}_{k}$. With the linear feedback controller $\xi_{k}=K_{k}\eta_{k}$, the inclusion $\eta_{k}\in\mathcal{E}_{Q_{k}}$ implies that $q_{k}$ is in a compact set $\mathcal{Q}$ that is given as 
\begin{align*}
\delta\mathcal{Q}_{k} & =\mathcal{E}_{C_{k}^{cl}Q_{k}(C_{k}^{cl})^{\top}}\oplus\{F_{k}w_{k}\mid\norm{w_{k}}_{2}\leq1\},\\
\mathcal{Q}_{k} & =\{\bar{q}_{k}\}\oplus\delta {\mathcal{Q}}_{k},\quad\forall\, k\in\mathcal{N}_{0}^{N-1},
\end{align*}
where $C_{k}^{cl}\coloneqq C+DK_{k}$. The assumption that the function $f$ is locally Lipschitz implies that the function ${\phi_{k}}$ is locally Lipschitz as well. Thus, for the compact (closed and bounded) set $\mathcal{Q}_{k}$, there exists a $\gamma_{k}$ such that
\begin{align*}
\norm{{\phi_{k}}(q_{k})-{\phi_{k}}(\bar{q}_{k})}_{2} & \leq\gamma_{k}\norm{q_{k}-\bar{q}_{k}}_{2}, \\
\forall\, q_{k}\in\mathcal{Q}_{k}, & \forall\, k\in\mathcal{N}_{0}^{N-1}.
\end{align*}
Considering them together, the closed-loop system becomes
\begin{subequations}\label{eq:diff_closedloop}
\begin{align}
\eta_{k+1} & =A_{k}^{cl}\eta_{k}+F_{k}w_{k}+E\delta p_{k}, \\
\delta q_{k} & =C_{k}^{cl}\eta_{k}+Gw_{k}, \\
\norm{\delta p_{k}}_{2} & \le\gamma_{k}\norm{\delta q_{k}}_{2}, \\
\norm{w_k} & \leq 1, \\
 & \delta q_{k}\in\delta\mathcal{Q},\quad\forall\, k\in\mathcal{N}_{0}^{N-1},
\end{align}
\end{subequations}
where $A_{k}^{cl}\coloneqq A_{k}+B_{k}K_{k}$. 

\subsubsection{Invariance of a quadratic funnel }
Consider a scalar-valued quadratic Lyapunov function $V$ defined by 
\begin{equation}
V(k,\eta)=\eta_{k}^{\top}Q_{k}^{-1}\eta_{k}.\label{eq:Lyapunov_function}
\end{equation}
For the closed-loop system model \eqref{eq:diff_closedloop}, we aim to design $\{Q_{k}\}_{k=0}^{N},\{K_{k}\}_{k=0}^{N-1}$ that satisfies the following quadratic stability condition:
\begin{subequations}\label{eq:Lyapunov_condition}
\begin{align}
V(k+1,\eta_{k+1}) & \leq\alpha V(k,\eta_{k}),\label{eq:Lyapunov_contraction}\\
\forall\,\norm{\delta p_{k}}_{2} & \le\gamma_{k}\norm{\delta q_{k}}_{2},\label{eq:Lyapunov_nonlinearity}\\
\forall\, V(k,\eta_{k}) & \geq\norm{w_{k}}^{2}_2,\label{eq:Lyapunov_disturbances}\\
 & \forall\, k\in\mathcal{N}_{0}^{N-1}\nonumber 
\end{align}
\end{subequations}
where $0<\alpha<1$. The above condition ensures the quadratic stability \eqref{eq:Lyapunov_contraction} whenever the locally Lipschitz property of the nonlinearity ${\phi_{k}}$ expressed in \eqref{eq:Lyapunov_nonlinearity} holds. The condition \eqref{eq:Lyapunov_disturbances} exists to obtain the invariance property of the funnel under the presence of the bounded disturbance $w_{k}$. In the rest of this subsection, we construct a condition that implies the stability condition \eqref{eq:Lyapunov_condition}. In the following corollary, we also show that the derived LMI condition ensures the invariance property of the funnel.

\begin{theorem}
Suppose that there exists $Q_{k}\in\mathbb{S}_{++}^{n_{x}}$, $Y_{k}\in\mathbb{R}^{n_{u}\times n_{x}}$,
$\nu_{k}^{p}>0$, $\lambda_{k}^{w}>0$, and $0<\alpha<1$ such
that $\lambda_k^w < \alpha$ and the following matrix inequality holds for all $k\in\mathcal{N}_{0}^{N-1}$:
\begin{equation}
\left[\begin{array}{ccccc}
\alpha Q_{k}-\lambda_{k}^{w}Q_{k} & * & * & * & *\\
0 & \nu_{k}^{p}I & * & * & *\\
0 & 0 & \lambda_{k}^{w}I & * & *\\
A_{k}Q_{k}+B_{k}Y_{k} & \nu_{k}^{p}E_{k} & F_{k} & Q_{k+1} & *\\
C_{k}Q_{k}+D_{k}Y_{k} & 0 & G_{k} & 0 & \nu_{k}^{p}\frac{1}{\gamma_{k}^{2}}I
\end{array}\right]\succeq0.\label{eq:LMI}
\end{equation}
Then the Lyapunov condition \eqref{eq:Lyapunov_condition} holds for the closed loop system \eqref{eq:diff_closedloop} with $K_{k}=Y_{k}Q_{k}^{-1}$.
\end{theorem}

\begin{proof}
With the closed-loop system \eqref{eq:diff_closedloop}, the condition
\eqref{eq:Lyapunov_condition}
holds if there exists a $\lambda_{k}^{p}>0$, $\lambda_{k}^{w}>0$,
and $0<\alpha\,{<}\,1$ such that the following inequality holds by $\mathcal{S}$-procedure \cite{yakubovich1997s}  for all $\eta_{k}\in\mathbb{R}^{n_{x}},w_{k}\in\mathbb{R}^{n_{w}},\delta p\in\mathbb{R}^{n_{p}}$:
\begin{align}
V(k+1,\eta_{k+1})-\alpha V(k,\eta_k)+\lambda_{k}^{w}(V(k,\eta_{k})-\norm{w_{k}}_{2}^{2}) +\lambda_{k}^{p}(\gamma_{k}^{2}\norm{\delta q_{k}}_{2}^{2}-\norm{\delta p_{k}}_{2}^{2})\leq0.\label{eq:attractive_in_V}
\end{align}
This is equivalent to
\begin{align*}
\left[\begin{array}{ccc}
A_{k}^{cl} & E_{k} & F_{k}\end{array}\right]^{\top}Q_{k+1}^{-1}\left[\begin{array}{ccc}
A_{k}^{cl} & E_{k} & F_{k}\end{array}\right]
+\lambda_{k}^{p}\left[\begin{array}{ccc}
C_{k}^{cl} & 0 & G_{k}\\
0 & I & 0
\end{array}\right]^{\top}\left[\begin{array}{cc}
\gamma_{k}^{2}I & 0\\
0 & -I
\end{array}\right]\left[\begin{array}{ccc}
C_{k}^{cl} & 0 & G_{k}\\
0 & I & 0
\end{array}\right]
-\left[\begin{array}{ccc}
\alpha Q_{k}^{-1} & * & *\\
0 & 0 & *\\
0 & 0 & 0
\end{array}\right]
+\lambda_{k}^{w}\left[\begin{array}{ccc}
Q_{k}^{-1} & * & *\\
0 & 0 & *\\
0 & 0 & -I
\end{array}\right]\preceq0.
\end{align*}
With the appropriate re-arrangement and applying Schur complement, we obtain
\begin{align*}
\left[\begin{array}{ccccc}
H_{k}^{1} & * & * & * & *\\
0 & \lambda_{k}^{p}I & * & * & *\\
0 & 0 & \lambda_{k}^{w}I & * & *\\
Q_{k+1}^{-1}A_{k}^{cl} & Q_{k+1}^{-1}E_{k} & Q_{k+1}^{-1}F_{k} & Q_{k+1}^{-1} & *\\
C_{k}^{cl} & 0 & G_{k} & 0 & H_{k}^{2}
\end{array}\right] & \succeq0
\end{align*}
where $H_{k}^{1}$ and $H_{k}^{2}$ are given by
\begin{align*}
H_{k}^{1} & =\alpha Q_{k}^{-1}-\lambda_{k}^{w}Q_{k}^{-1},\\
H_{k}^{2} & =(\lambda_{k}^{p})^{-1}\frac{1}{\gamma_{k}^{2}}I.
\end{align*}
Multiplying both sides by $\text{diag}\{Q_{k},\lambda_{p}^{-1}I,Q_{k+1},I\}$ yields
\begin{align*}
\left[\begin{array}{ccccc}
\alpha Q_{k}-\lambda_{k}^{w}Q_{k} & * & * & * & *\\
0 & \nu_{k}^{p}I & * & * & *\\
0 & 0 & \lambda_{k}^{w}I & * & *\\
A_{k}^{cl}Q_{k} & \nu_{k}^{p}E_{k} & F_{k} & Q_{k+1} & *\\
C_{k}^{cl}Q_{k} & 0 & G_{k} & 0 & \nu_{k}^{p}\frac{1}{\gamma_{k}^{2}}I
\end{array}\right]\succeq0,
\end{align*}
where $\nu_{k}^{p}=(\lambda_{k}^{p})^{-1}$. Finally, expanding $A_{k}^{cl}$ and $C_{k}^{cl}$ completes the proof.
\end{proof}
\begin{corollary}
The condition \eqref{eq:LMI} in Theorem 1 implies the following invariance condition for all $k\in\mathcal{N}_{0}^{N-1}$:
\begin{subequations}\label{eq:cor_invariance}
\begin{align}
V(k+1,\eta_{k+1}) & \leq1,\\
\forall V(k,\eta_{k}) & \leq1,\label{eq:invariance_V}\\
\forall\norm{\delta p_{k}}_2 & \leq\gamma_{k}\norm{\delta q_{k}}_2,\\
\forall\norm{w_{k}}_{2} & \leq1.\label{eq:invariance_disturbance}
\end{align}
\end{subequations}
\end{corollary}
\begin{proof}
Observe that \eqref{eq:attractive_in_V} can be equivalently written as
\begin{align*}
V(k+1,\eta_{k+1})-\alpha+(\alpha-\lambda_{k}^{w})(1-V(k,\eta_{k})) +\lambda_{k}^{w}(1-\norm{w_{k}}_{2}^{2})+\lambda_{k}^{p}(\gamma_{k}^{2}\norm{\delta q_{k}}_{2}^{2}-\norm{\delta p_{k}}_{2}^{2})\leq0.
\end{align*}
This implies $V(k+1,\eta_{k+1}) \leq\alpha$ with \eqref{eq:invariance_V}-\eqref{eq:invariance_disturbance}
since $0<\lambda_{k}^{w}<\alpha$ and $\lambda_{k}^{p}>0$. This is
sufficient for the invariance condition \eqref{eq:cor_invariance} since $\alpha<1$.
\end{proof}
\noindent Notice that the matrix inequality \eqref{eq:LMI} is a LMI once $\alpha$
and $\lambda_{k}^{w}$ are fixed.

\subsubsection{Computing the funnel via SDP}
The goal of computing the CIF is to bound the effects of disturbances going forward in time by minimizing the size of the funnel while satisfying the invariance and the feasibility of the boundary conditions. To this end, the funnel computation is posed as the following SDP:
\begin{subequations}\label{eq:funnel_update}
\begin{align}
\underset{
\begin{array}{c} \scriptscriptstyle Q_{k}, {\mu^Q_{k}},\forall\, k\in\mathcal{N}_{0}^{N}, \\[-0.1cm]
\scriptscriptstyle Y_{k},{\mu^K_{k}},\nu_{k}^{p},\forall\, k\in\mathcal{N}_{0}^{N-1}
\end{array}}{\operatorname{minimize}} & \change{w_Q\sum_{k=0}^{N} \mu^Q_{k} + w_K\sum_{k=0}^{N-1} \mu^K_{k} + \sum_{k=0}^{N-1}J_{trf}(Q_{k},Y_{k}) }\label{eq:funnelopt_cost}\\[-0.1cm]
\operatorname{subject~to}~~& 
Q_{k}\preceq{\mu^Q_{k}}I,\forall\, k\in\mathcal{N}_{0}^{N},\label{eq:funnelopt_minQ}\\
 & \left[\begin{array}{cc}
{\mu^K_{k}}I & Y_{k}\\
Y_{k}^{\top} & Q_{k}
\end{array}\right]\succeq0,\forall\, k\in\mathcal{N}_{0}^{N-1},\label{eq:funnelopt_minK}\\
 & \operatorname{condition~}\eqref{eq:LMI},\\
&\change{\left[\begin{array}{cc}
( (b_{i}^{x})_{k} - (a_{i}^{x})_{k}^{\top}\bar{x}_{k} )^2  & (a_{i}^{x})_{k}^{\top} Q_{k}^\top\\
Q_{k} (a_{i}^{x})_{k} & Q_{k}
\end{array}\right]\succeq0,i=1,\ldots,m_{x},}\label{eq:funnel_state_const} \\
&\change{\left[\begin{array}{cc}
( (b_{j}^{u})_{k} - (a_{j}^{u})_{k}^{\top}\bar{u}_{k} )^2  & (a_{j}^{u})_{k}^{\top} Y_{k}^\top\\
Y_{k} (a_{j}^{u})_{k} & Q_{k}
\end{array}\right]\succeq0,j=1,\ldots,m_{u},}\label{eq:funnel_input_const}\\
 & Q_{0}\succeq Q_{i},\quad Q_{N}\preceq Q_{f},\label{eq:funnelopt_initial}
\end{align}
\end{subequations}
where \eqref{eq:funnelopt_minK} is equivalent to \eqref{eq:main_minimizeK} which can be derived by Schur complement with $Y_{k}=K_{k}Q_{k}$. The LMI constraints in \eqref{eq:funnel_state_const}-\eqref{eq:funnel_input_const} are the funnel feasibility conditions that are equivalent to \eqref{eq:lin_feasibility} that can be derived by Schur complement. The cost $J_{trf}$ is given as
\begin{align*}
J_{trf}=w_{trf}\sum_{k=0}^{N-1}\left(\norm{Q_{k}-\hat{Q}_{k}}_{F}^2+\norm{Y_{k}-\hat{Y}_{k}}_{F}^2\right),
\end{align*}
where $w_{trf}\in\mathbb{R}$ is a user-defined parameter, $\norm{\cdot}_{F}$ is the Frobenius norm, and $\hat{Y}_{k}=\hat{K}_{k}\hat{Q}_{k}$ for all $k\in\mathcal{N}_{0}^{N-1}$. This cost, similar to the trust region penalty $J_{tr}$, penalizes the difference between the current solution $\{Q_{k},Y_{k}\}_{k=0}^{N-1}$ and the previous solution $\{\hat{Q}_{k},\hat{Y}_{k}\}_{k=0}^{N-1}$ which is beneficial for the better convergence performance.

\change{The choice of parameters in the proposed method affects the performance of the control law in (3). The weights $w_Q$ and $w_K$ in (7a) balances the size of the state funnel $\mathcal{E}_{Q_k}$ and input funnel $\mathcal{E}_{K_k Q_k K_k^\top}$. For example, a relatively larger $w_Q$ compared to $w_K$ drives the algorithm to put more effort on minimizing the size of the state funnel over the input funnel, and vice versa. The choices of the decay rate $\alpha$ and the slack variable $\lambda_k^w$ resulted from $\mathcal{S}$-procedure in (18) also affects the control performance. As the decay rate decreases, the controller places greater emphasis on faster convergence to the nominal trajectory due to the condition outlined in (17a). Likewise, the larger $\lambda_k^w$ places more emphasis on the convergence to the nominal trajectory. This is attributed to the term $V(k+1,\eta_{k+1})-\alpha V(k,\eta_k)$ in (19) becoming smaller (more negative)  as $\lambda_{k}^{w}(V(k,\eta_{k})-\norm{w_{k}}_{2}^{2})$ increases.
}

\subsubsection{Local Lipschitz constant estimation via sampling}
To compute the LMI \eqref{eq:LMI}, the Lipschitz constant $\gamma_{k}$ in \eqref{eq:Lyapunov_condition} should be available. We estimate the Lipschitz constant by employing a sampling method. It is worth mentioning that the sampling method for the estimation of the Lipschitz constant $\gamma_{k}$ brings about an algebraic loop: to estimate the Lipschitz constant $\gamma_{k}$, the funnel variables $Q_{k}$ and $K_{k}$ should be available, whereas the computation $Q_{k}$ and $K_{k}$ in \eqref{eq:LMI} requires the constant $\gamma_{k}$. However, a well-behaved iterative scheme with the sampling method for $\gamma_{k}$ can make the funnel computation converge \cite{reynolds2021funnel}.

By sampling a set of $N_{s}$ pairs of state and disturbance $\{\eta_{k}^{s},w_{k}^{s}\}_{s=1}^{N_{s}}$ from the ellipsoid $\mathcal{E}_{Q}$ and the set $\{w\in\mathbb{R}^{n_{w}}\mid\norm w_{2}\leq1\}$, respectively, we compute
\begin{align}
\delta_{k}^{s} & =\frac{\norm{p_{k}^{s}-\bar{p}_{k}}}{\norm{q_{k}^{s}-\bar{q}_{k}}},\quad {s}=1,\ldots,N_{s},\label{eq:direct_evaluation_gamma}
\end{align}
where $p_{k}^{s}$ and $q_{k}^{s}$ are computed by \eqref{eq:converted_model}. Depending on the discretization method,  only $Ep_{k}$ might be available instead of $p_k$. So, it might not be possible to compute \eqref{eq:direct_evaluation_gamma}. In that case, we instead solve the following optimization to obtain the value $\delta_{k}^{s}$:
\begin{subequations}\label{eq:indirect_evaluation_gamma}
\begin{align}
\delta_{k}^{s}=
\underset{\scriptscriptstyle \Delta}{\operatorname{minimize}}~~ & \norm{\Delta}_{2} \\
\operatorname{subject~to} \ & \eta_{k+1}^{s}-A_{k}^{cl}\eta_{k}^{s}-F_{k}w_{k}^{s}+\bar{x}_{k+1}  - f(\bar{x}_{k},\bar{u}_{k},0) = E\Delta(C_{k}^{cl}\eta_{k}^{s}+G_{k}{w_{k}^{s}}),
\end{align}
\end{subequations}
where $\Delta\in\mathbb{R}^{n_{p}\times n_{q}}$.  After obtaining $\delta_{k}^{s}$ by \eqref{eq:direct_evaluation_gamma} or \eqref{eq:indirect_evaluation_gamma}, the following maximization operation is performed to estimate the local Lipschitz constant:
\begin{align}
\gamma_{k}=
\underset{\scriptscriptstyle s=1,\ldots,N_{s}}{\operatorname{maximize}}~
\delta_{k}^{s},\quad\forall\, k\in\mathcal{N}_{0}^{N-1}. \label{eq:gamma_update}
\end{align}
It is worth noting that the disadvantage of the illustrated sampling-based method is that the computed $\gamma_k$ might be lower than the true level of nonlinearity. To handle this issue, one may able to use a probabilistic approach for overestimating the local Lipschitz constant from samples provided in \citenum{chou2021model}.

\change{Another way to estimate the local Lipschitz constant is to use the optimization-based approach provided in Section 6.5.1 of \citenum{reynolds2020computational}. To illustrate, we consider
\begin{subequations}\label{eq:outer_optimization_lipschitz}
\begin{align}
\Gamma_{k}^*= 
\underset{\scriptscriptstyle \eta_k, w_k}{\operatorname{maximize}}~~ & \frac{1}{2} \delta_k^*(\eta_k,w_k)^2 \\
\operatorname{subject~to} \ & \eta_k Q_k^{-1} \eta_k \leq 1,  \\
& \norm{w_k}_2 \leq 1, 
\end{align}
\end{subequations}
where
\begin{subequations}\label{eq:inner_optimization_lipschitz} 
\begin{align}
\delta_{k}^*(\eta_k,w_k)=
\underset{\scriptscriptstyle \Delta}{\operatorname{minimize}}~~ & \norm{\Delta}_{2} \\
\operatorname{subject~to} \ & \eta_{k+1}-A_{k}^{cl}\eta_{k} - F_{k}w_{k}+\bar{x}_{k+1} - f(\bar{x}_{k},\bar{u}_{k},0) = E\Delta(C_{k}^{cl}\eta_{k} + G_{k}{w_{k}}). \label{eq:inner_const}
\end{align}
\end{subequations}
The inner optimization \eqref{eq:inner_optimization_lipschitz} aims to find the smallest matrix $\Delta$ in terms of the matrix 2-norm for the given $\eta_k, w_k$, and the outer optimization \eqref{eq:outer_optimization_lipschitz} finds the values of $\eta_k, w_k$ that maximize $\delta_k^*$. After solving these optimization problems for each $k$, the local Lipschitz constant can be obtained by computing
\begin{align*}
    \gamma_k = \sqrt{2 \Gamma_k^*}.
\end{align*}
To make the outer optimization computationally tractable, one could potentially utilize an analytic upper bound for the problem's optimal value. The constraint \eqref{eq:inner_const} can be rewritten as
\begin{subequations}
    \begin{align}
        y(\eta_k,w_k) &= E\Delta(C_{k}^{cl}\eta_{k} + G_{k}{w_{k}}), \\
        &= \underbrace{\left[\begin{array}{ccc} E(e_1^\top(C_{k}^{cl}\eta_{k} + G_{k}{w_{k}})) & \cdots &
        E(e_{n_q}^\top(C_{k}^{cl}\eta_{k} + G_{k}{w_{k}}))
        \end{array}\right]}_{\coloneqq H(\eta_k,w_k)} \Vec{\Delta}_k,
    \end{align}
\end{subequations}
where  $y(\eta_k,w_k) = \eta_{k+1}-A_{k}^{cl}\eta_{k} - F_{k}w_{k}+\bar{x}_{k+1} - f(\bar{x}_{k},\bar{u}_{k},0)$ and $\Vec{\Delta}_k \in \mathbb{R}^{n_p n_q}$ is a concatenated vector that stacks the columns of $\Delta_k$. Then, consider the following optimization:
\begin{align}
    \delta_k^* (\eta_k,w_k) = \underset{\scriptscriptstyle \Vec{\Delta}_k}{\operatorname{minimize}}
    \ \norm{\Vec{\Delta}_k}_2 \quad \operatorname{subject~to} \ \ y(\eta_k,w_k)=H(\eta_k,w_k)\Vec{\Delta}_k. \label{eq:approximate_inner_optimization}
\end{align}
Since \eqref{eq:approximate_inner_optimization} is a minimum-norm least squares problem, the solution is $\Vec{\Delta}_k = H^\dagger (\eta_k,w_k)y(\eta_k,w_k)$ \cite{boyd2004convex}  where $H^\dagger$ represents the peusoinverse of matrix $H$. The optimal value $\delta_k^*$ of \eqref{eq:approximate_inner_optimization} is the upper bound of the inner optimization \eqref{eq:inner_optimization_lipschitz} as $\norm{\Delta}_2\leq \norm{\Delta}_F = \norm{\Vec{\Delta}_k}_2$ where $\norm{\cdot}_F$ is the Frobenius norm. Then, the outer optimization \eqref{eq:outer_optimization_lipschitz} can be transformed into
\begin{subequations} \label{eq:approximate_outer_optimization}
\begin{align}
    \Gamma_{k}^*= 
    \underset{\scriptscriptstyle \eta_k, w_k}{\operatorname{maximize}}~~ & \frac{1}{2} y(\eta_k,w_k)^\top H^\dagger(\eta_k,w_k)^\top H^\dagger(\eta_k,w_k) y(\eta_k,w_k)  \\
    \operatorname{subject~to} \ & \eta_k Q_k^{-1} \eta_k \leq 1,  \\
    & \norm{w_k}_2 \leq 1, 
\end{align}
\end{subequations}
More details in the derivation and the computation results can be found Section 6.5.1 in \citenum{reynolds2020computational}.}

\subsection{\change{Algorithm details and summary}} \label{subsec:algorithm}

To start the algorithm, we need to generate an initial guess hat is used as a reference trajectory for the first iteration. It is worth noting that the initial guess does not need to be feasible to constraints for the proposed method. The first way is to employ a straight-linear interpolation \cite{malyuta2022convex} for the initial nominal trajectory $\{x_{k}\}_{k=0}^{N},\{u_{k}\}_{k=0}^{N-1}$. Then the feedback gain $\{K_{k}\}_{k=0}^{N-1}$ can be obtained by solving a discrete-time linear quadratic regulator problem with a linearized model of \eqref{eq:nonlinear} evaluated around the nominal trajectory. The initial guess for the ellipsoid variable $\{Q_{k}\}_{k=0}^{N}$ can then be set to a diagonal matrix having user-defined diameters. \change{The second way to generate the initial guess is to use the separate synthesis; the nominal trajectory is generated by the SCP algorithm without considering the funnel. This provides the dynamically-feasible trajectory that can be used as the initial guess of the nominal trajectory $\{{x}_{k}\}_{k=0}^{N},\{{u}_{k}\}_{k=0}^{N-1}$ for the proposed method. Then, the feedback gain $\{K_{k}\}_{k=0}^{N-1}$ and the ellipsoid variable $\{Q_{k}\}_{k=0}^{N}$ can be obtained via solving \eqref{eq:funnel_update} with $w_{trf}=0$ while ignoring the funnel feasibility. The second way is more systematical since it exploits the result of the separate synthesis and hence gives a better initial guess compared to the solution computed by the straight-line interpolation in the first way.}

To set the stopping criteria, we define $\Delta_{vc}, \Delta_{dyn}, \Delta_{T}$ and $\Delta_{F}$ as
\begin{align*}
\Delta_{vc} &= \sum_{k=0}^{N-1}\norm{v_{k}}_{1}, \ 
\Delta_{dyn} = \sum_{k=0}^{N-1} \norm{f(t_k,\bar{x}_{k},\bar{u}_{k},0)-\bar{x}_{k+1}}_2, \\
\Delta_{T} & =\norm{x_{N}-\hat{x}_{N}}_{2}^{2}+\sum_{k=0}^{N-1}\norm{x_{k}-\hat{x}_{k}}_{2}^{2}+\norm{u_{k}-\hat{u}_{k}}_{2}^{2}, \\
\Delta_{F} & =\norm{Q_{N}-\hat{Q}_{N}}_{F}^{2}+\sum_{k=0}^{N-1}\norm{Q_{k}-\hat{Q}_{k}}_{F}^{2}+\norm{Y_{k}-\hat{Y}_{k}}_{F}^{2}.
\end{align*}
Then the stopping criteria is given as the following logical statement:
\begin{align}
(\Delta_{vc}<\Delta_{vc}^{tol})\wedge
(\Delta_{dyn}<\Delta_{dyn}^{tol})\wedge
(\Delta_{T}<\Delta_{T}^{tol})\wedge
(\Delta_{F}<\Delta_{F}^{tol}), \label{eq:stopping_criterion}
\end{align}
where $\Delta_{vc}^{tol}, \Delta_{dyn}^{tol}, \Delta_{T}^{tol}$ and $\Delta_{F}^{tol}$ are user-defined tolerance parameters. The proposed algorithm is summarized in Algorithm \ref{alg:alg1}.

\begin{algorithm}
\begin{algorithmic}
\Require{($\hat{x}_k,\hat{u}_k\,\hat{Q}_k,\hat{K}_k$)}
\For{$i=1\ldots N_{max}$}
    \State{optimize $\bar{x}_k,\bar{u}_k$ by \eqref{eq:traj_update}}
	\State{estimate $\gamma_k$ via \eqref{eq:gamma_update} or \eqref{eq:approximate_outer_optimization}}
    \State{optimize $Q_k,K_k$ by \eqref{eq:funnel_update}}
    \If{\eqref{eq:stopping_criterion} is True}
        \State{break}
    \EndIf
    \State{update $(\hat{x}_k,\hat{u}_k\,\hat{Q}_k,\hat{K}_k)\leftarrow (\bar{x}_k,\bar{u}_k,Q_k,K_k)$}
\EndFor
\Ensure{$(\bar{x}_k,\bar{u}_k,Q_k,K_k)$}
\end{algorithmic}
\caption{Joint synthesis}
\label{alg:alg1}
\end{algorithm}

\change{While the convergence guarantee of the proposed method has not been a focus of this paper, one can construct a safety alternative that is assured not to diverge by modifying the proposed algorithm with results from \citenum{cartis2011evaluation, reynolds2021funnel}. Instead of updating the trajectory and the funnel sequentially in each iteration, the safety approach performs updating only the trajectory with a fixed funnel until convergence of the nominal trajectory is achieved. This part of the process, being solely trajectory optimization, benefits from the established convergence results in \citenum{cartis2011evaluation}. The subsequent phase involves computing the Lipschitz constant and updating the funnel with the computed nominal trajectory, with convergence analysis from Theorem 6.12 of \citenum{reynolds2021funnel}. Since each phase of the safety approach has a guaranteed convergence, it prevents the overall solution from diverging.}

\section{Numerical simulation}
\label{sec4}

\change{In this section, we validate the proposed method via two robotic applications with a unicycle model and a 6-DoF free-flying spacecraft.} For both examples, we used an Apple MacBook Pro having M1 Pro with 8-core CPU, and the simulation result can be reproduced using the code available at {{\href{https://github.com/taewankim1/joint\_synthesis}{\texttt{https://github.com/taewankim1/joint\_synthesis}}}}.

\subsection{Unicycle model}
\change{We consider the motion of a unicycle-type model  under different disturbance conditions, represented by $w_1$ and $w_2$
\begin{subequations}\label{eq:unicycle_model}
\begin{align}
\text{Model $\textrm{I}$: }\left[\begin{array}{c}
\dot{r_x}\\
\dot{r_y}\\
\dot{\theta}
\end{array}\right] & =\left[\begin{array}{c}
u_{v}\cos\theta+0.1w_{1}\\
u_{v}\sin\theta+0.1w_{2}\\
u_{\theta}
\end{array}\right], \quad\text{Model $\textrm{II}$: }\left[\begin{array}{c}
\dot{r_x}\\
\dot{r_y}\\
\dot{\theta}
\end{array}\right] =\left[\begin{array}{c}
(u_{v}+0.1w_{1})\cos\theta\\
(u_{v}+0.1w_{1})\sin\theta\\
u_{\theta} + 0.1w_{2}
\end{array}\right], \\
\text{Model $\textrm{III}$: }
\left[\begin{array}{c}
\dot{r_x}\\
\dot{r_y}\\
\dot{\theta}
\end{array}\right] & =\left[\begin{array}{c}
u_{v}\cos(\theta+0.03w_{1})\\
u_{v}\sin(\theta+0.03w_{1})\\
u_{\theta} + 0.05w_{2}
\end{array}\right],
\end{align}
\end{subequations}
where $r_x$, $r_y$, and $\theta$ are $x$-axis position, $y$-axis position, are heading angle, respectively, and $u_{v}\in\mathbb{R}$ is velocity and $u_{\theta}\in\mathbb{R}$ is angular velocity. The scalars $w_1$ and $w_2$ represent disturbances or model mismatch. Model $\textrm{I}$ depicts direct disturbances on the translational motion, Model $\textrm{II}$ introduces disturbances affecting both the velocity and rotational control inputs, and Model $\textrm{III}$ captures disturbances influencing the orientation and the rotation control. These models are considered to have a comprehensive understanding of how the system behaves according to different types of the disturbances. It is worth noting that in Model $\textrm{II}$ and Model $\textrm{III}$, the disturbances introduce additional nonlinearities to the system, while in Model $\textrm{I}$, they appear as linear additive terms.}

\begin{figure}

\begin{center}
    \centering
     \begin{subfigure}{0.49\textwidth}
     \centering
     \includegraphics[width=1.0\textwidth,trim={6cm 0cm 6cm 1cm},clip]{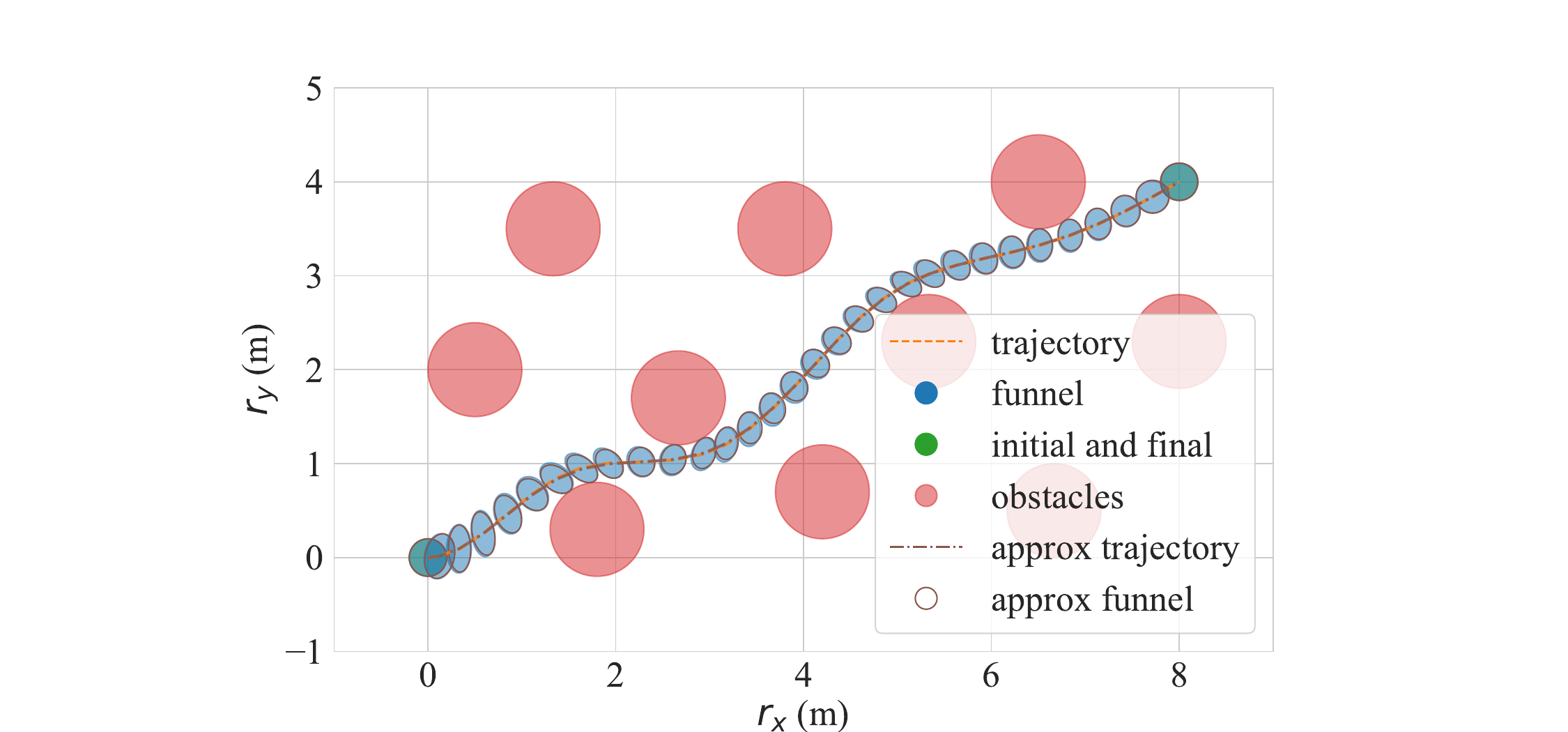}
     \end{subfigure}
     \begin{subfigure}{0.49\textwidth}
     \centering
     \includegraphics[width=1.0\textwidth,trim={6cm 0cm 6cm 1cm},clip]{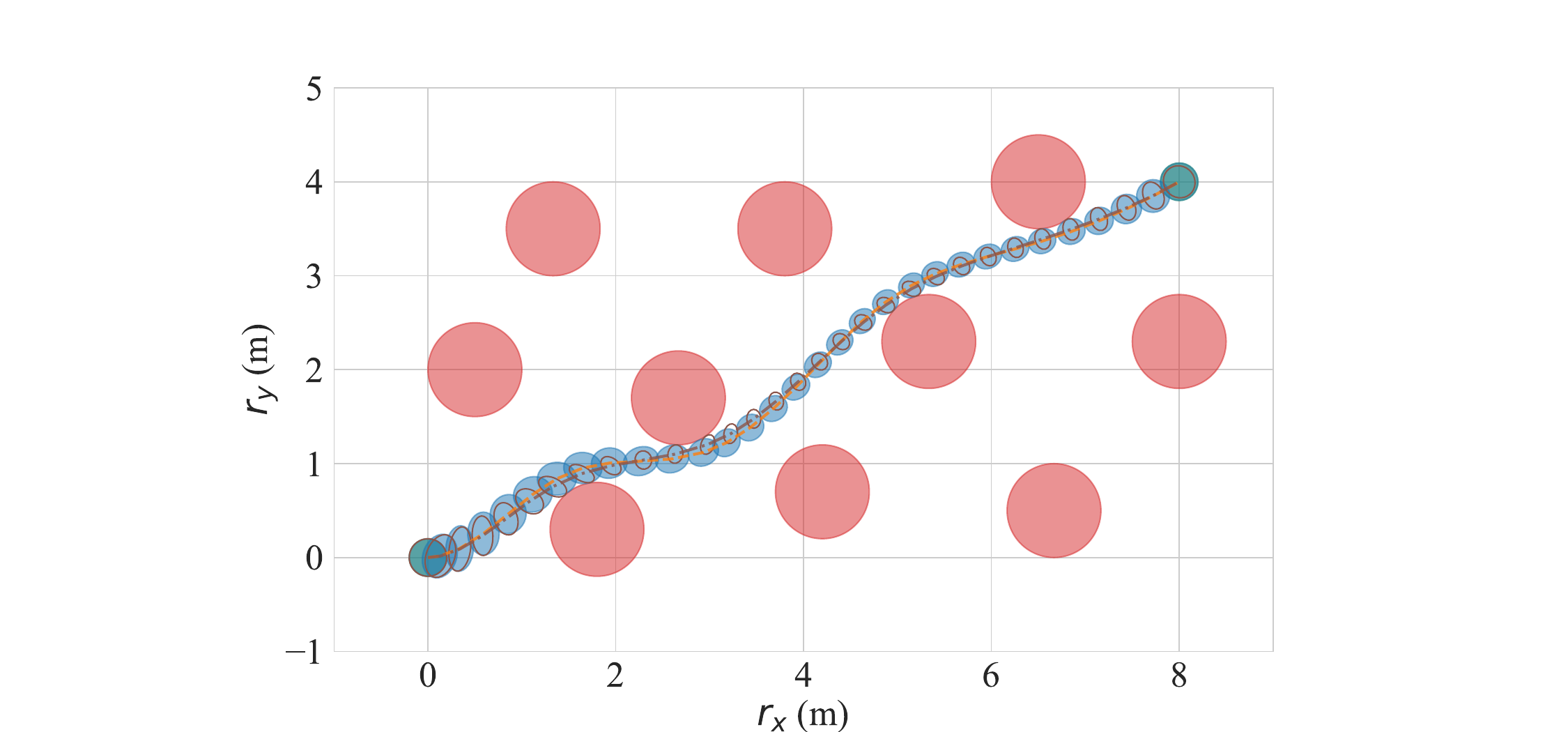}
     \end{subfigure}
     \begin{subfigure}{0.49\textwidth}
     \centering
     \includegraphics[width=1.0\textwidth,trim={6cm 0cm 6cm 1cm},clip]{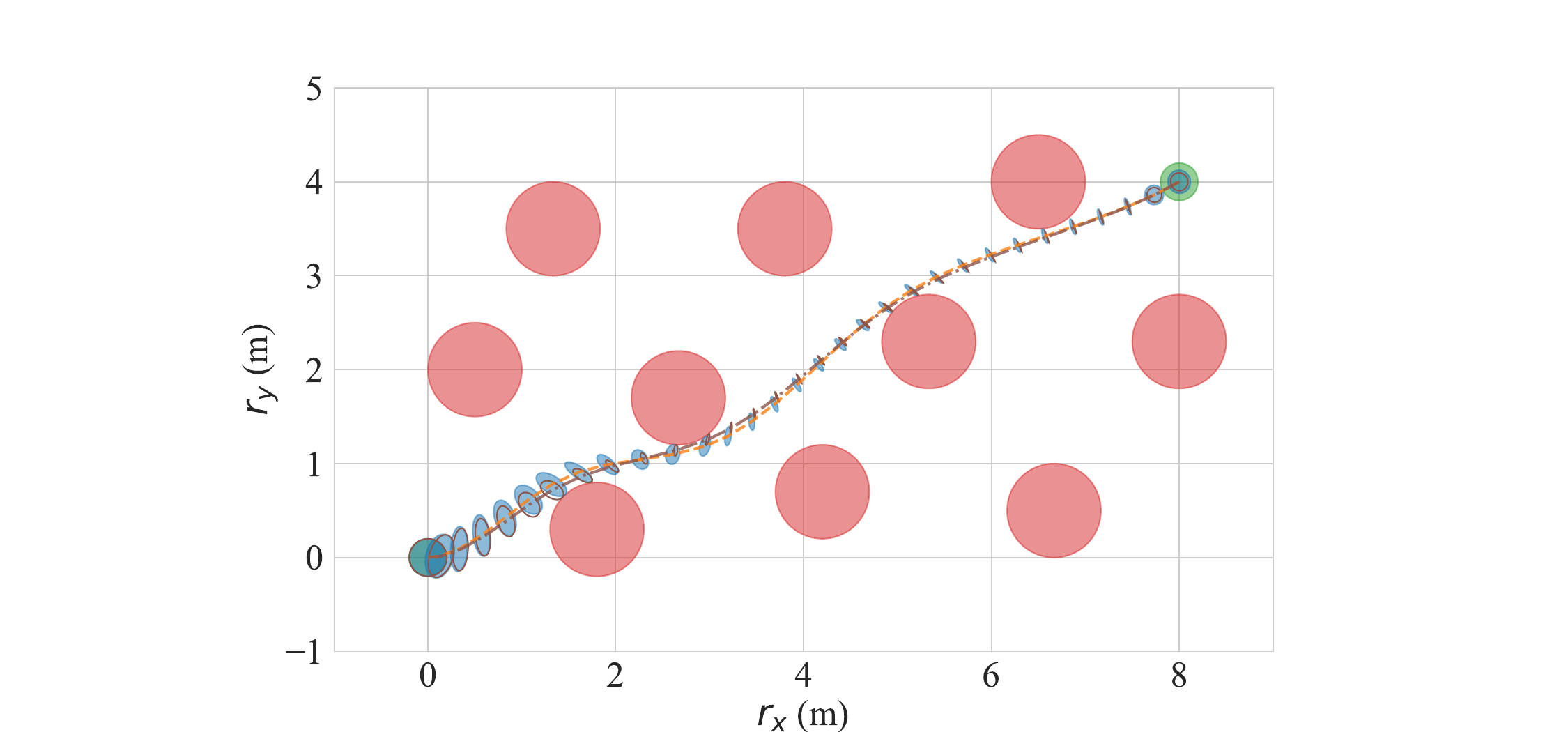}
     \end{subfigure}
\caption{
\change{Nominal trajectories and synthesized funnels (projected on position coordinates) of Model $\textrm{I}$ (top-left), Model  $\textrm{II}$ (top-right), and Model  $\textrm{III}$ (bottom). Each figure shows the nominal trajectory (orange line), the projection of the state ellipsoid in the funnel (blue ellipse), and the approximated funnel generated with the linear closed-loop system (brown ellipse).}
} 
\label{fig:traj_result}
\end{center}
\end{figure}

For all unicycle models, we consider $N=30$ nodes evenly distributed over a time horizon of \change{$10$ s i.e., $t_{0}=0$ and $t_{f}=10$}. The continuous-time model \eqref{eq:unicycle_model} is discretized by following a variational approach \cite[Chap. 10.4]{diehl2011numerical} to obtain the matrices $A_k, B_k, F_k$. The initial boundary set $\mathcal{X}_0$ and the final boundary set $\mathcal{X}_f$ in \eqref{eq:boundary_details} have the following parameters: $x_0 = [0,0,0]^\top$, \change{$Q_i = Q_f = \text{diag}([0.2^2$ (m), $0.2^2$ (m), $10^2$ (deg)$]^\top)$}, and $x_f = [8,4,0]^\top$. There are multiple circular obstacles the unicycle robot should avoid, which leads to nonconvex constraints on the state represented in set $\mathcal{X}$. All obstacles have a diameter of 1.0m, and their center positions are illustrated in Fig.~\ref{fig:traj_result}. \change{The input constraints for the set $\mathcal{U}$ are given as: $0\leq u_v \leq 1.5$ and $\abs{u_\theta} \leq 1.0$ (rad).} The cost function for the trajectory $J_t$ is a quadratic function of the input given by $u_v^2 + u_\theta^2$. \change{Both weight parameters $w_Q, w_K$ in \eqref{eq:maincost} are chosen as 1.} The decay rate $\alpha$ is set as $0.99$ and the parameter $\lambda^w_k$ is set as $0.2$ for all $k$. The tolerance parameters $\Delta_{vc}^{tol}, \Delta_{dyn}^{tol}, \Delta_T^{tol}$ and $\Delta_F^{tol}$ are all $10^{-8}$. The number of samples $N_s$ used for the Lipschitz constant $\gamma_k$ estimation is set as 100, for each $k$, so a total of 3,000 samples are used for each iteration. \change{We use an interior-point method solver, Clarabel, for both the trajectory update \eqref{eq:traj_update} and the funnel update \eqref{eq:funnel_update}, using CVXPY in Python.}

\begin{figure}
\begin{center}
    \centering
     \begin{subfigure}{0.3\textwidth}
     \centering
     \includegraphics[width=1.0\textwidth,trim={0cm 0cm 0cm 0cm},clip]{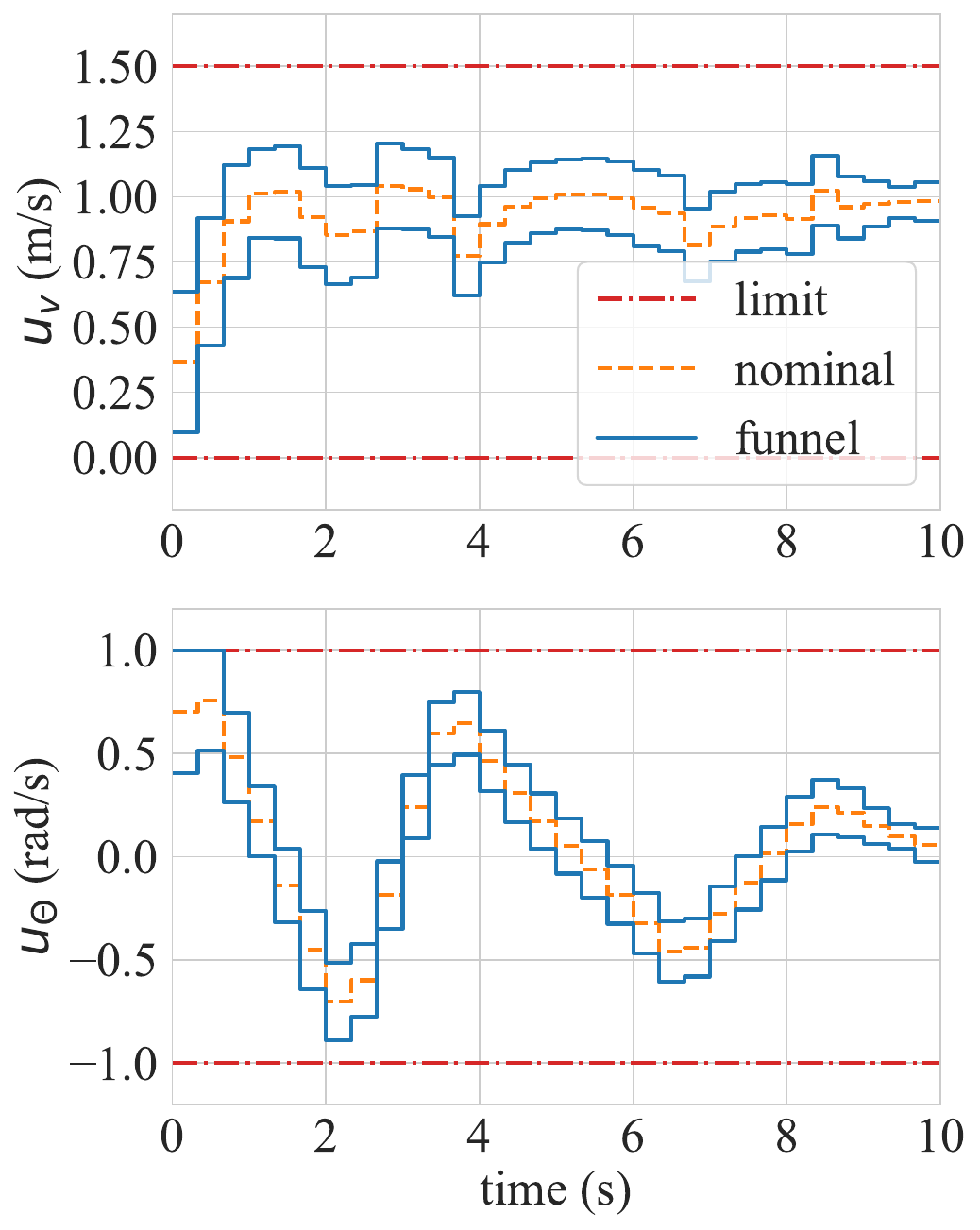}
     \end{subfigure}
     \begin{subfigure}{0.3\textwidth}
     \centering
     \includegraphics[width=1.0\textwidth,trim={0cm 0cm 0cm 0cm},clip]{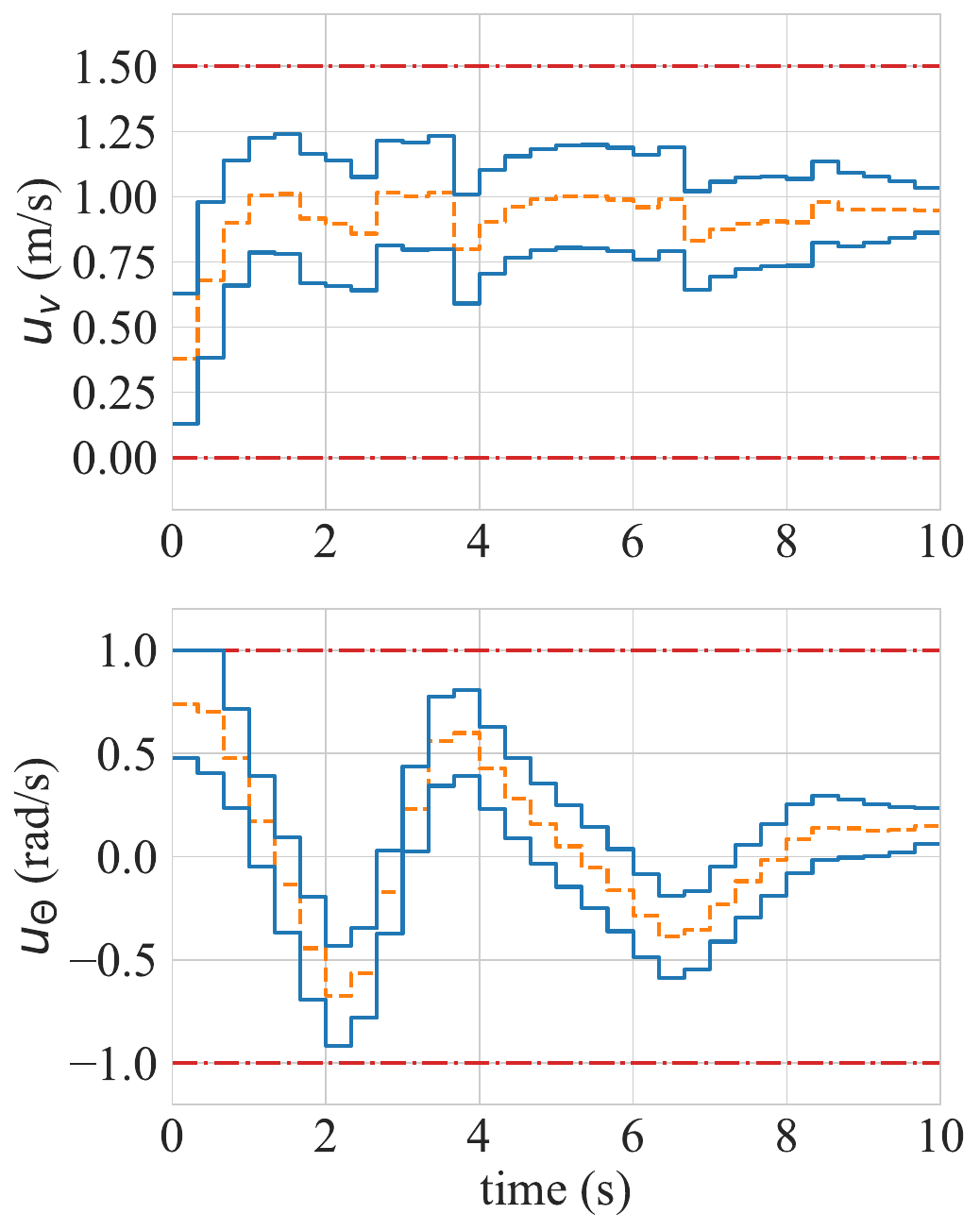}
     \end{subfigure}
     \begin{subfigure}{0.3\textwidth}
     \centering
     \includegraphics[width=1.0\textwidth,trim={0cm 0cm 0cm 0cm},clip]{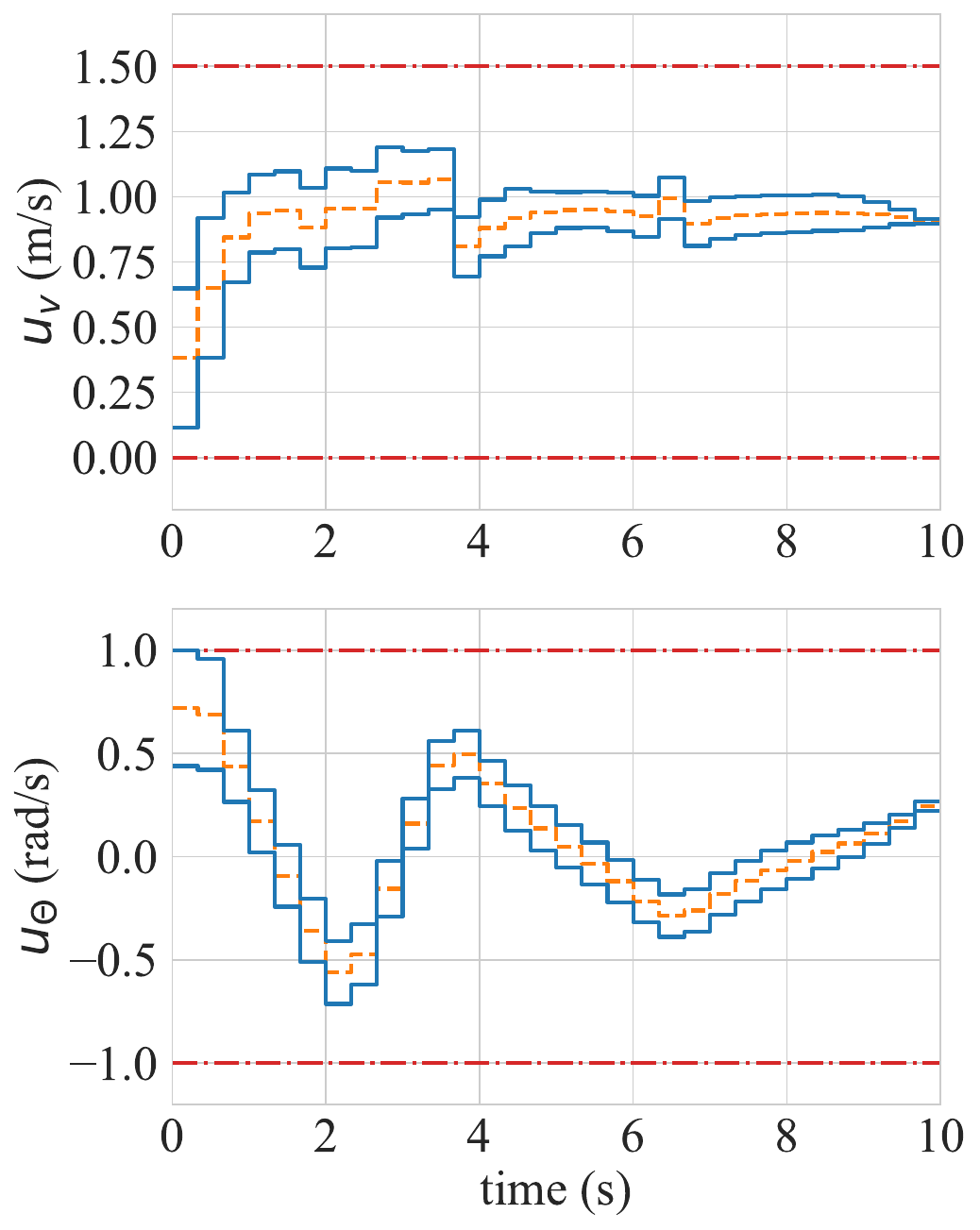}
     \end{subfigure}
\caption{
\change{Nominal trajectories and synthesized input funnels (projected on each input coordinate) of Model $\textrm{I}$ (left), Model  $\textrm{II}$ (middle), and Model  $\textrm{III}$ (right). The zeroth-order hold on the input is used to generate the nominal trajectory.}
} 
\label{fig:input_results}
\end{center}
\end{figure}

\begin{figure}
\begin{center}
\includegraphics[width=1.0\textwidth,trim={0cm 0cm 0cm 0cm},clip]{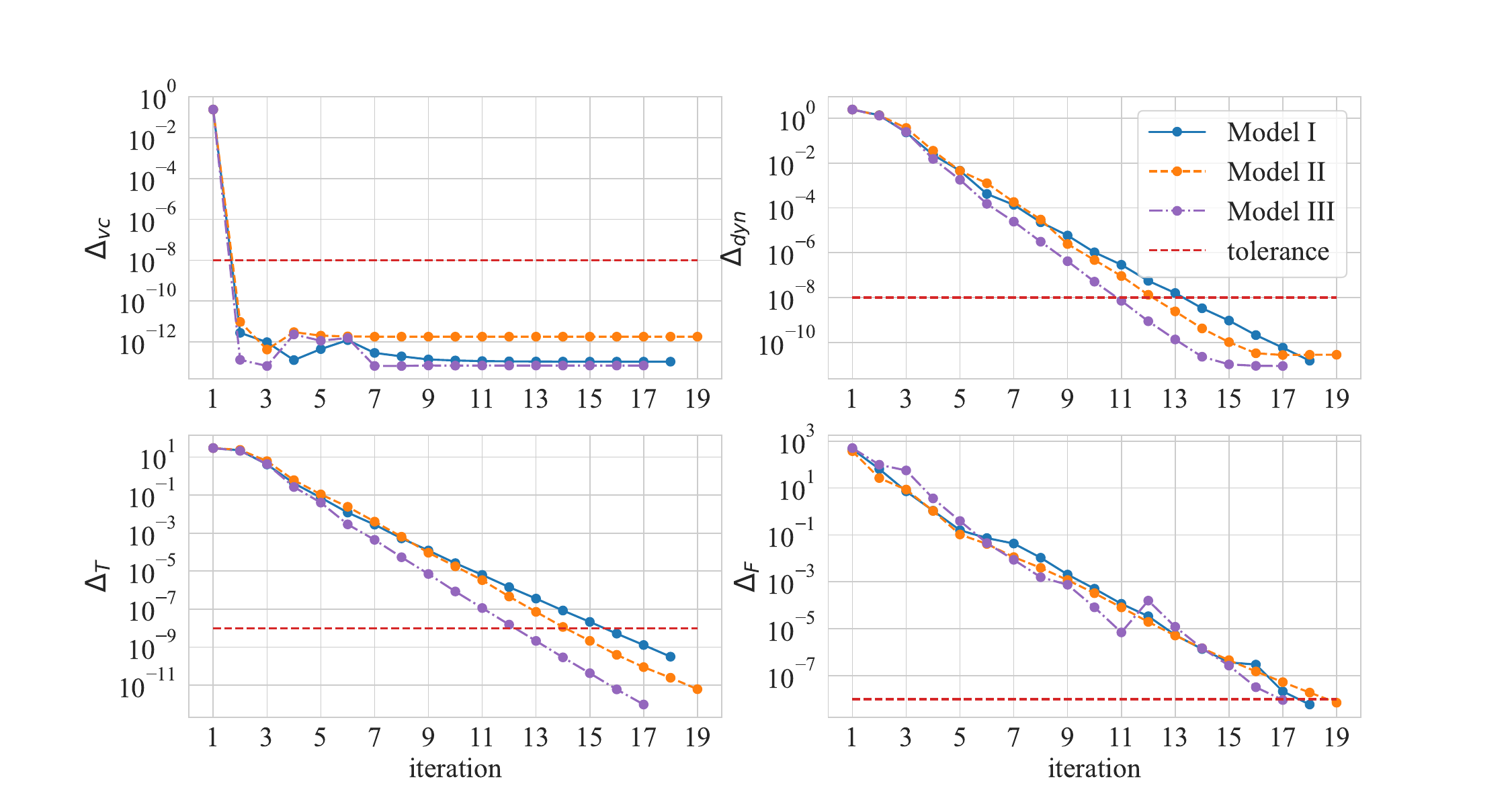}  
\caption{
\change{Convergence performance of the proposed method for the unicycle models.}
} 
\label{fig:convergence_results}
\end{center}
\end{figure}

\begin{figure}
\begin{center}
    \centering
     \begin{subfigure}{0.32\textwidth}
     \centering
     \includegraphics[width=1.0\textwidth,trim={0cm 0cm 0cm 0cm},clip]{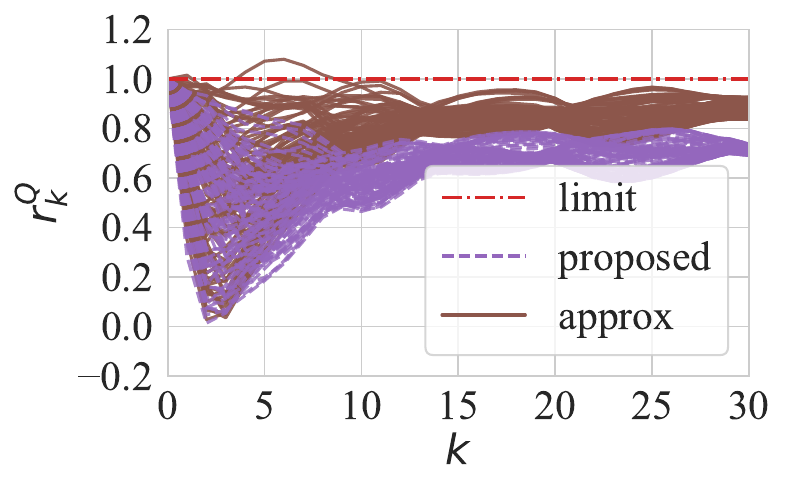}
     \end{subfigure}
     \begin{subfigure}{0.32\textwidth}
     \centering
     \includegraphics[width=1.0\textwidth,trim={0cm 0cm 0cm 0cm},clip]{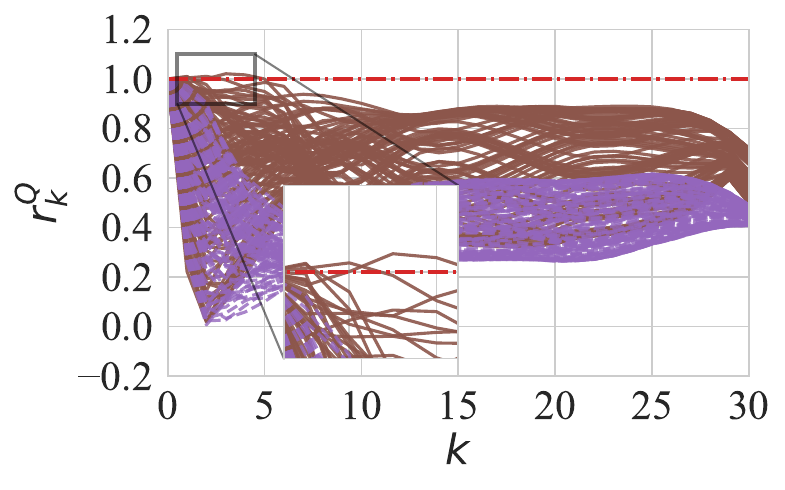}
     \end{subfigure}
     \begin{subfigure}{0.32\textwidth}
     \centering
     \includegraphics[width=1.0\textwidth,trim={0cm 0cm 0cm 0cm},clip]{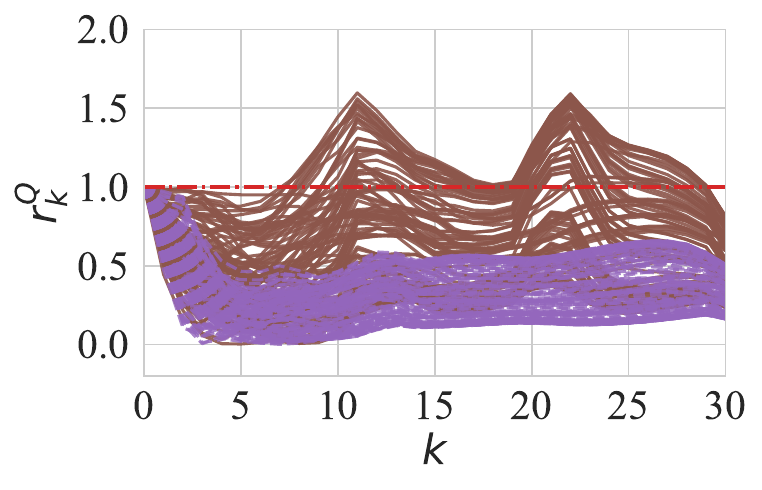}
     \end{subfigure}
\caption{
\change{Invariance property tests for Model $\textrm{I}$ (left), Model  $\textrm{II}$ (middle), and Model  $\textrm{III}$ (right).}
} 
\label{fig:comparison}
\end{center}
\end{figure}

To test the invariance and the feasibility properties, we sample 100 points at the surface of the ellipse $\mathcal{E}_{Q_k}$ at $k=0$, and then generate the corresponding 100 trajectories with the nonlinear dynamics \eqref{eq:nonlinear} and the control law \eqref{eq:control_law}  under the presence of the disturbances. In this generation process, we randomly set the disturbance $w=(w_1,w_2)$ such that $\norm{w}_2=1$ and keep them constant during the entire horizon for each sample. 
Note that making the disturbance constant during the entire horizon increases the impact of the disturbance compared to varying the disturbance randomly for each interval. The computed nominal trajectory and the CIF \change{for all unicycle models \eqref{eq:unicycle_model}} are depicted in Fig.~\ref{fig:traj_result}, and the input results are given in Fig.~\ref{fig:input_results}. The test results of the invariance property for the trajectory samples is given in Fig.~\ref{fig:comparison} where the radius $r^Q_k$ is defined as 
\begin{align}
    r^Q_k\coloneqq(x_k^s-\bar{x}_k)^\top Q^{-1}_k (x_k^s-\bar{x}_k) \label{eq:r_k}
\end{align}
for each sample $s$ and time $k$. The result shows that the nominal trajectory and the CIF satisfy the invariance and feasibility conditions. \change{For the initial guess, we employ the first method illustrated in Sec.~\ref{subsec:algorithm} using the straight-linear interpolation.} The convergence performance in Fig.~\ref{fig:convergence_results} shows that the proposed approach makes the trajectory and the CIF satisfy the tolerances as the iteration count increases. Table~\ref{tb:computation_time} summarizes the average computational time of each subproblem within the iterations.

To obtain a baseline solution to compare against, we compute an approximate funnel that is generated with the linear closed-loop system where the higher-order terms are ignored, that can be established by setting $E_k=0$ and $\gamma_k=0$ in \eqref{eq:diff_closedloop}, as considered in \citenum{messerer2021efficient,manchester2019robust}. It is worth noting that the approximate funnel, which is used for the comparison with the proposed method, can yield more optimal solutions compared to \citenum{messerer2021efficient,manchester2019robust} under the linear approximation. This is because the approximate funnel is computed by simultaneously optimizing the linear feedback gains and the invariant set parameters as decision variables, whereas \citenum{messerer2021efficient} determines the invariant set variables by the uncertainty forward equation and \citenum{manchester2019robust} sets the feedback gains by solving a discrete-time linear quadratic regulator problem. The approximate nominal trajectory and funnel are depicted in Fig.~\ref{fig:traj_result}, and the invariance test $r^Q_k \leq 1$ in \eqref{eq:r_k} with the trajectory samples is given in Fig.~\ref{fig:comparison}. We can see that the value of $r^Q_k$ for the approximate funnel is greater than 1 \change{especially for Model \textrm{III} since the bounded disturbances contribute the nonlinearity of the system.} This violation shows that the approximate CIF does not necessarily guarantee the invariance property for the original nonlinear system, which can result in safety issues for safety-critical nonlinear systems. As the contribution of higher order terms increase, e.g., for large Lipshitz constants, these violations can become more pronounced.

\begin{table}
\caption{Average computational time (s) for each iteration}
\label{tb:computation_time}
\centering{}%
\begin{tabular}{cccc}
\toprule 
Subproblem &  Trajectory update & Estimate $\gamma_k$ & Funnel update	\tabularnewline
\midrule
Model \textrm{I} & 0.026 & 0.863 & 0.577 \\ 
\midrule 
Model \textrm{II} & 0.026 & 0.895 & 0.691 \\  
\midrule 
Model \textrm{II} & 0.025 & 0.904 & 0.880 \\  
\bottomrule
\end{tabular}
\end{table}

\subsection{6-DoF free-flying spacecraft}

\change{We consider the following 6-DoF free-flying spacecraft dynamics \cite{malyuta2022convex,lew2020chance} under the presence of disturbances:
\begin{subequations}
    \label{eq:freeflyer_dyn}
    \begin{align}
        \dot{r}_{\mathcal{I}} &= v_{\mathcal{I}}, \\
        \dot{v}_{\mathcal{I}} &= m^{-1} (T_{\mathcal{I}} + \beta_T w_T), \\
        \dot{\Phi} &= R(\Phi)\omega_{\mathcal{B}}, \\
        \dot{\omega}_{\mathcal{B}} &= J^{-1}(M_{\mathcal{B}} + \beta_M w_M - \omega_{\mathcal{B}} \times J\omega_{\mathcal{B}}).
    \end{align}
\end{subequations}
The state of the system \eqref{eq:freeflyer_dyn} consists of the inertial position $r_{\mathcal{I}}\in \mathbb{R}^3$, the inertial velocity $v_{\mathcal{I}}\in \mathbb{R}^3$, the ZYX Euler angles $\Phi\in\mathbb{R}^3$, and the body angular velocity ${\omega}_{\mathcal{B}} \in \mathbb{R}^3$. The control input consists of the inertial thrust $T_{\mathcal{I}}\in\mathbb{R}^3$ and the body torque $M_{\mathcal{B}}\in\mathbb{R}^3$. The constant $m\in\mathbb{R}$ is the mass, matrix $J\in\mathbb{R}^{3\times3}$ is the inertia matrix, and $R$ maps the Euler angles to the rotation matrix
\begin{align*}
    R(\Phi) = \left[\begin{array}{ccc}
1 & \sin \phi \tan \theta & \cos\phi \tan \theta\\
0 & \cos \phi & -\sin \phi \\
0 & \sin \phi \sec \theta & \cos \phi \sec \theta 
\end{array}\right],
\end{align*}
with $\Phi = [\phi,\theta,\psi]^\top$. For the mass and inertia matrix parameters, we choose $m = 7.2$ (kg) and $J=0.1083 I_{3\times3}$ (kgm$^2$) where $I_{3\times 3}$ is the 3 by 3 identity matrix. The vector $w_T\in\mathbb{R}^3$ and $w_M \in \mathbb{R}^3$ affect the system as disturbances by acting on the control inputs $T_{\mathcal{I}}$ and $M_{\mathcal{B}}$ with coefficients $\beta_T=10^{-3}$ and $\beta_M=10^{-6}$. Further details about the free-flying system dynamics can be found in \citenum{malyuta2022convex}.}

\begin{figure}
    \centering
        \centering
         \begin{subfigure}[b]{0.4\textwidth}
         \centering
         \includegraphics[width=\textwidth,clip]{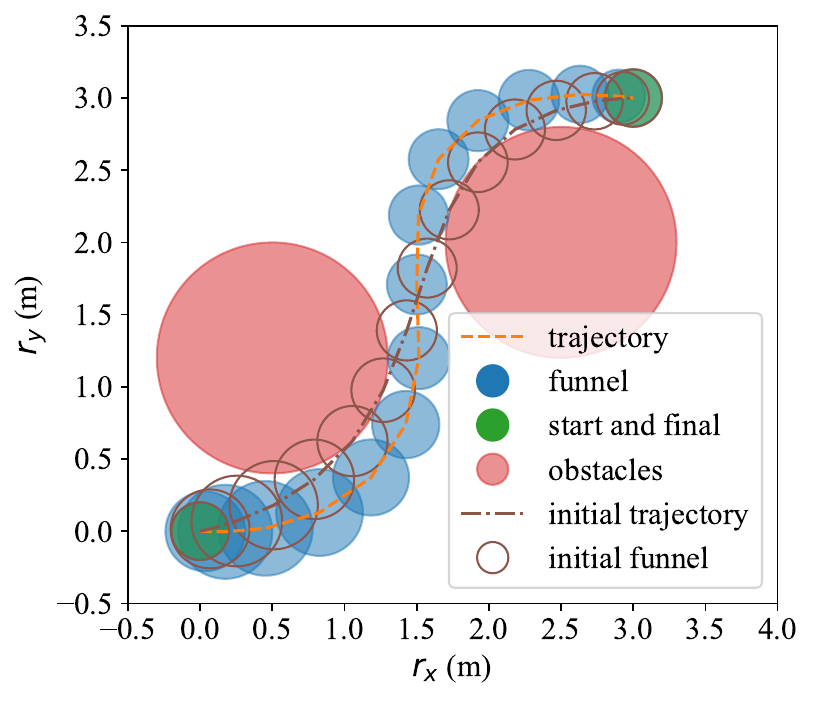}
         \end{subfigure}
         \begin{subfigure}[b]{0.5\textwidth}
         \centering
         \includegraphics[width=0.8\textwidth,trim={2cm 2cm 1cm 3cm},clip]{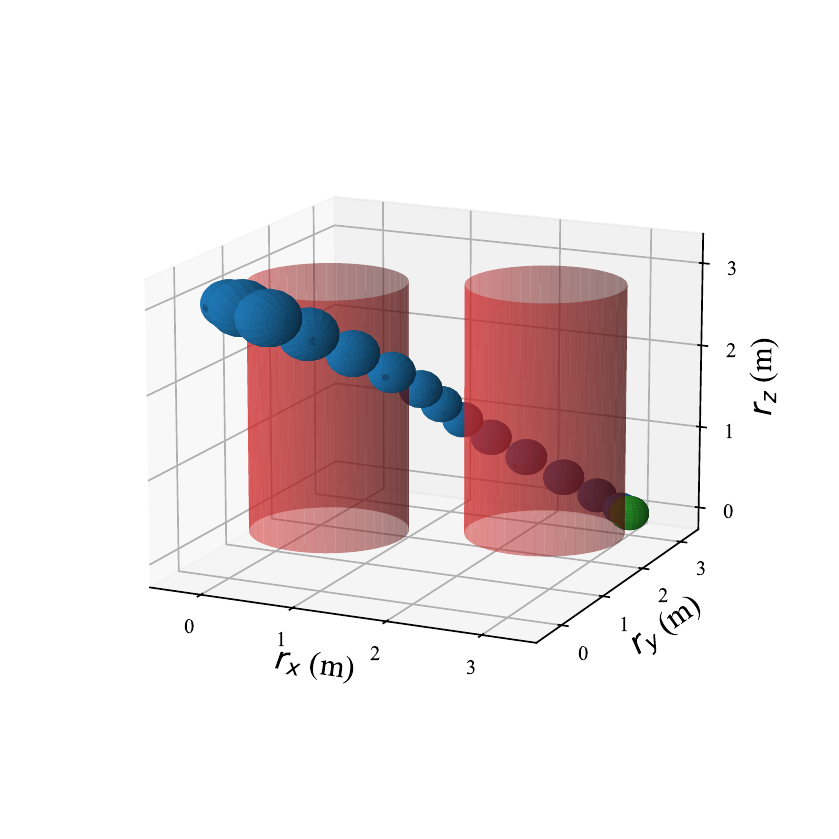}
         \end{subfigure}
\caption{Nominal trajectories and synthesized funnel projected on $r_xr_y$ positions (left) and $r_xr_yr_z$ positions (right), respectively, for the free-flying spacecraft.}
\label{fig:freeflyer_traj_result}
\end{figure}

\change{For the free-flying spacecraft example, we consider $N=15$ nodes evenly distributed over a time horizon of $200$ s. 
The initial boundary set $\mathcal{X}_0$ and the final boundary set $\mathcal{X}_f$ in \eqref{eq:boundary_details} have the following parameters: 
\begin{align*}
    x_0 &= [r_0^\top,v_0^\top,\Phi_0^\top,\omega_0^\top]^\top, r_0 = [0,0,3]^\top\text{(m)}, v_0=[0,0,0]^\top\text{(m/s)},
    \Phi_0 = [-30,25,5]^\top\text{(deg)},
    \omega_0 = [0,0,0]^\top\text{(deg/s)},\\
    x_f &= [r_f^\top,v_f^\top,\Phi_f^\top,\omega_f^\top]^\top,
    r_f = [3,3,0]^\top\text{(m)}, v_f=\Phi_f=\omega_f=[0,0,0]^\top, \\
    Q_i &= Q_f = \text{diag}\left(\left[0.2^2, 0.2^2, 0.2^2, 0.02^2,0.02^2,0.02^2,
    \left(\frac{5\pi}{180}\right)^2,\left(\frac{5\pi}{180}\right)^2,\left(\frac{5\pi}{180}\right)^2,
    \left(\frac{0.1\pi}{180}\right)^2,\left(\frac{0.1\pi}{180}\right)^2,\left(\frac{0.1\pi}{180}\right)^2\right]\right).
\end{align*}
As state constraints for the set $\mathcal{X}$, there are two cylindrical obstacles to avoid, and all obstacles have a diameter of 0.8m, and their center positions are illustrated in Fig.~\ref{fig:freeflyer_traj_result}. In addition, we have constraints on the velocity and the angular velocity that are $\norm{v_{\mathcal{I}}}_2 \leq 0.4$ (m/s) and $\norm{\omega_{\mathcal{B}}}_2 \leq 1$ (deg). The input constraints for the set $\mathcal{U}$ are given as: $ \norm{T_{\mathcal{I}}}_2 \leq 10$ (mN) and $ \norm{M_{\mathcal{B}}}_2 \leq 50$ ($\mu$Nm).
The decay rate $\alpha$ is set as $0.99$ and the parameter $\lambda^w_k$ is set as $0.1$ for all $k$. Similar to the unicycle examples, the tolerance parameters $\Delta_{vc}^{tol}, \Delta_{dyn}^{tol}, \Delta_T^{tol}$ and $\Delta_F^{tol}$ are all $10^{-8}$. The number of samples $N_s$ used for the Lipschitz constant $\gamma_k$ estimation is set as 256 for each $k$. We provide the initial guess using the second method illustrated in Sec.~\ref{subsec:algorithm} that uses the result of the separate synthesis, and these used initial trajectory and funnel are illustrated in Fig.~\ref{fig:freeflyer_traj_result}. Starting from the initial guess, the proposed algorithm converges at 6 iterations. We use Clarabel for the trajectory update \eqref{eq:traj_update} and MOSEK for the funnel update \eqref{eq:funnel_update}, using CVXPY in Python3. Mosek's solve time is observed to scale better than Clarabel's for large problem sizes. So, for the free-flyer system, we have used Mosek for the funnel update. The average computational time (s) of the trajectory update, the estimation of $\gamma_k$, and the funnel update at each iteration are 0.024, 2.698, and 10.072, respectively. }

\begin{figure}
\begin{center}
\includegraphics[width=0.8\textwidth,trim={0cm 0cm 0cm 0cm},clip]{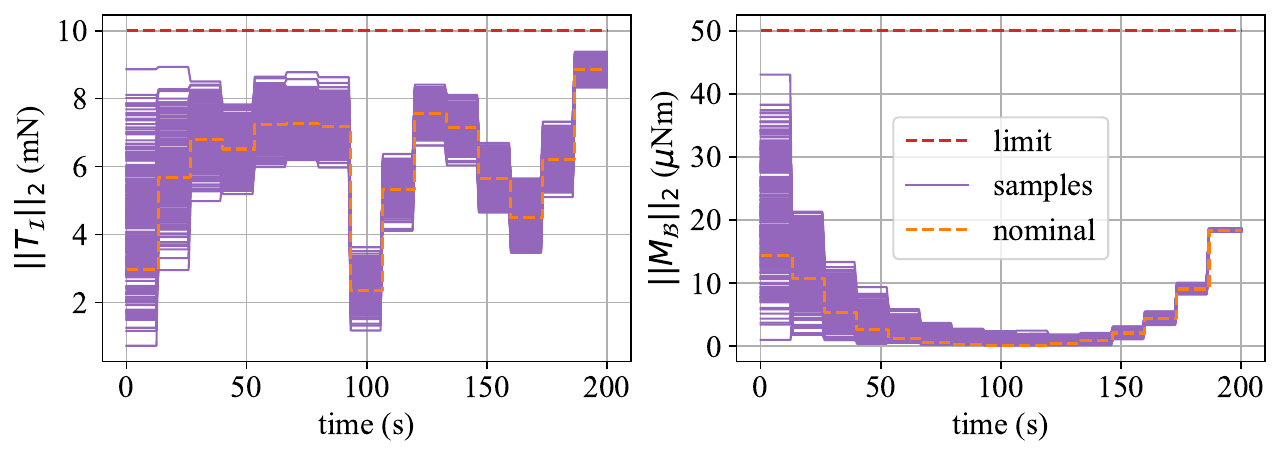}  
\caption{
The thrust and moment results of nominal trajectory and trajectory samples for the free-flying spacecraft.
} 
\label{fig:freeflyer_input_result}
\end{center}
\end{figure}

\begin{figure}
\begin{center}
\includegraphics[width=0.5\textwidth,trim={0cm 0cm 0cm 0cm},clip]{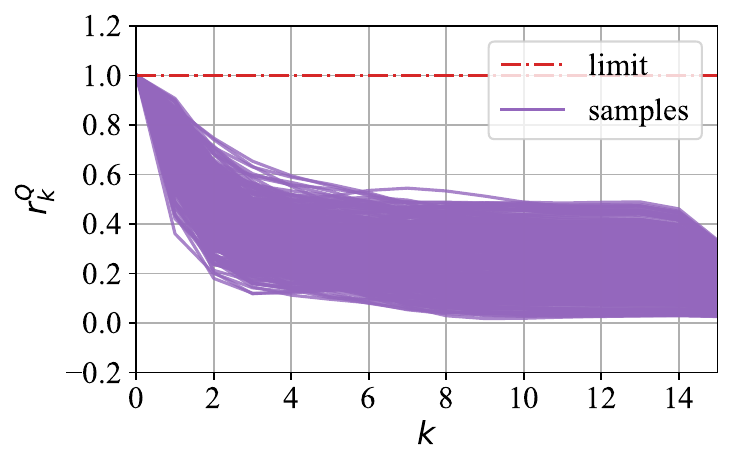} 
\caption{
Invariance property tests for trajectory samples of the free-flying spacecraft.
} 
\label{fig:freeflyer_invariance}
\end{center}
\end{figure}

\change{The results of the synthesized trajectory and funnel projected on position coordinates are illustrated in Fig.~\ref{fig:freeflyer_traj_result}. It is clear that the resulting funnel is feasible to the obstacle avoidance constraints although the initial guess has infeasible trajectory and funnel. Similar to the test performed for the unicycle models, to test the invariance of the synthesized trajectory and funnel, we sample 300 at the surface of the initial funnel $\mathcal{E}_{Q_0}$ and generate the corresponding 300 trajectories by propagating the system dynamics \eqref{eq:freeflyer_dyn} with the control law \eqref{eq:control_law} consisting of the open-loop input in the nominal trajectory and the feedback control from the funnel. For each sample, we randomly set the disturbance $w=[w_T^\top,w_M^\top]^\top$ such that $\norm{w}_2=1$. The input results of the nominal trajectory and the samples are illustrated in Fig.~\ref{fig:freeflyer_input_result}. We can see that the input history of the samples remain feasible within the given input constraints. Finally, we obtain the values of $r^Q_k$ by computing \eqref{eq:r_k} for each trajectory sample, illustrated in Fig.~\ref{fig:freeflyer_invariance}. The result shows that $r^Q_k$ for each sample and for all time $k$ maintain less than 1 and hence the trajectory samples remain inside the funnel.}

\section{Conclusions}
\label{sec5}
This paper introduces a method for joint trajectory optimization and funnel synthesis for locally Lipschitz nonlinear systems under the presence of disturbances. The proposed method has a recursive approach in which both nominal trajectory and funnel are iteratively updated. The trajectory update step optimizes the nominal trajectory to satisfy the feasibility of the funnel. Then, the funnel update step solves an SDP to guarantee the invariance property of the funnel. \change{The numerical evaluation for a unicycle model and a 6-DoF free-flying spacecraft shows that the converged trajectory and funnel satisfy the invariance and feasibility properties under the disturbances.}

\section{Acknowledgments}
We thank Dr. Taylor P. Reynolds and Prof. Mehran Mesbahi for their helpful inputs on funnel synthesis algorithms. This work is supported in part by Boeing under Grant 2021-PD-PA-471, Office of Navel Research under Grant N00014-20-1-2288, and AFOSR under Grant FA9550-20-1-0053. Taewan Kim is supported by UW+Amazon Science Hub fellowship.


\bibliography{wileyNJD-AMA}%




\end{document}